\documentclass{article}
\usepackage{latexsym,amsfonts,amsmath,amsthm,amssymb,makeidx}
\usepackage[title]{appendix}
\usepackage{CJK,CJKnumb,CJKulem,times,dsfont,ifthen,mathrsfs,latexsym,amsfonts, color}
\usepackage{amsmath,amsthm,makeidx,fontenc,amssymb,bm,graphicx,psfrag,listings, curves,extarrows}

\let\oldbibliography\thebibliography
\renewcommand{\thebibliography}[1]{%
\oldbibliography{#1}%
\setlength{\itemsep}{0pt}%
}

\makeindex
\newtheorem{definition}{Definition}[section]
\newtheorem{theorem}{Theorem}[section]
\newtheorem{lemma}{Lemma}[section]
\newtheorem{corollary}{Corollary}[section]

\newtheorem{remark}{Remark}[section]

\newcommand{\s}{\section}

\newcommand{\la}{\lambda}

\newcommand{\pa}{\partial}

\newcommand{\R}{\mathbb R}
\newcommand{\al}{\alpha}

\newcommand{\bt}{\begin{theorem}}
\newcommand{\et}{\end{theorem}}
\newcommand{\bl}{\begin{lemma}}
\newcommand{\el}{\end{lemma}}
\newcommand{\bd}{\begin{definition}}
\newcommand{\ed}{\end{definition}}
\newcommand{\bc}{\begin{corollary}}
\newcommand{\ec}{\end{corollary}}
\newcommand{\bp}{\begin{proof}}
\newcommand{\ep}{\end{proof}}
\newcommand{\bx}{\begin{example}}
\newcommand{\ex}{\end{example}}
\newcommand{\bi}{\begin{exercise}}
\newcommand{\ei}{\end{exercise}}
\newcommand{\bo}{\begin{prop}}
\newcommand{\eo}{\end{prop}}
\newcommand{\br}{\begin{remark}}
\newcommand{\er}{\end{remark}}
\newcommand{\be}{\begin{equation}}
\newcommand{\ee}{\end{equation}}
\newcommand{\ba}{\begin{align}}
\newcommand{\ea}{\end{align}}
\newcommand{\bn}{\begin{enumerate}}
\newcommand{\en}{\end{enumerate}}
\newcommand{\bg}{\begin{align*}}
\newcommand{\bcs}{\begin{cases}}
\newcommand{\ecs}{\end{cases}}

\newcommand{\bean}{\begin{eqnarray*}}
\newcommand{\eean}{\end{eqnarray*}}

\numberwithin{equation}{section}

\begin{document}

\title{\bf Classification of the stable solution to the  fractional $(2<s<3)$  Lane-Emden equation\thanks{Partially supported by NSFC of China and  NSERC of Canada. E-mails: luosp14@mails.tsinghua.edu.cn(Luo);\;  jcwei@math.ubc.ca (Wei);\; wzou@math.tsinghua.edu.cn(Zou)}}
\date{}
\author{\\{\bf  Senping Luo$^{1}$,\;\;  Juncheng Wei$^{2}$ \;\; and \; Wenming Zou$^{3}$}\\
\footnotesize {\it  $^{1,3}$Department of Mathematical Sciences, Tsinghua University, Beijing 100084, China}\\
\footnotesize {\it  $^{2}$Department of Mathematics, University of British Columbia, Vancouver, BC V6T 1Z2,
Canada}
}

\maketitle
\begin{center}
\begin{minipage}{120mm}
\begin{center}{\bf Abstract}\end{center}

We classify the stable solutions (positive or sign-changing, radial or not) to  the following nonlocal Lane-Emden equation:
\noindent
\begin{equation}\nonumber
(-\Delta)^s u=|u|^{p-1}u\;\;\;\;\hbox{in}\;\;\;\;\; \R^n
\end{equation}
for $2<s<3$.


\vskip0.10in


\end{minipage}
\end{center}

\vskip0.10in
\section{Introduction and Main results}

Consider  the stable solution of the following equation
 \begin{equation}\label{8LE}
(-\Delta)^s u=|u|^{p-1}u\;\;\;\;\hbox{in}\;\;\;\;\;  \R^n,
\end{equation}
where $(-\Delta)^s$ is the fractional Laplacian operator for $2<s<3$.

 \vskip0.1in
 The motivation of  studying    such an equation is originated from the classical   Lane-Emden equation
 \begin{equation}\label{LZ=000123}
 -\Delta  u=|u|^{p-1}u \;\;\;\;\hbox{in}\;\;\;\;\;  \R^n
\end{equation}
and its parabolic counterpart, which  have played a crucial   role in the development of nonlinear PDEs in the last
decades. These  arise  in astrophysics and Riemannian geometry.  The pioneering  works on  Eq.\eqref{LZ=000123} were  contributed by   R. Fowler \cite{Fowler=1,Fowler=2}. Later, the  ground-breaking  result on equation \eqref{LZ=000123} is the  fundamental    Liouville-type theorems  established by  Gidas and Spruck \cite{Gidas-Sp=1981}, they claimed   that the Eq. \eqref{LZ=000123} has no positive solution whenever $p\in (1, 2^\ast-1)$, where $2^\ast=2n/(n-2)$ if $n\geq 3$ and
 $2^\ast=\infty$ if $n\leq 2.$  The critical case  $p=2^\ast-1$,     Eq.\eqref{LZ=000123}  has a unique positive solution  up to translation  and rescaling  which is radial and explicitly formulated, see Caffarelli-Gidas-Spruck \cite{Caffarelli1989}. Since then
many experts in partial differential equations devote to the above equations for various
parameters $s$ and $p$.

 \vskip0.1in

For the nonlocal case of $0 < s < 1,$  a counterpart of the classification results of Gidas and Spruck \cite{Gidas-Sp=1981}, and Caffarelli-
Gidas-Spruck  \cite{Caffarelli1989}  holds for the fractional Lane-Emden equation  \eqref{8LE}, see the works due to   Li \cite{YLi04} and Chen-Li-Ou \cite{Chen-Li-Ou}. In these cases, the Sobolev exponent is given by $P_S(n, s)=(n+2s)/(n-2s)$ if $n>2s$, and otherwise $P_S(n, s)=\infty$.

\vskip0.123in
Recently, for the nonlocal case of $0<s<1$, Davila, Dupaigne and Wei in \cite{Wei0=1} gave a complete classification of finite Morse index solution of \eqref{8LE}; for the nonlocal case of $1<s<2$, Fazly and Wei in \cite{Wei1=2} gave a complete classification of finite Morse index solution of \eqref{8LE}.
For the local cases $s=1$ and $s=2$, such kind of  classification is proved by Farina in \cite{Farina2007} and Davila, Dupaigne, Wang and Wei in \cite{Wei=2}, respectively. For the case $s=3$, the Joseph-Lundgren exponent (for the triharmonic Lane-Emden equation) is obtained
and classification is proved by in \cite{LWZ2016=3}.

 \vskip0.1in

However, when $2<s<3$,  the equation \eqref{8LE} has not been  considered  so far. In this paper we classify the stable solution of \eqref{8LE}.
\vskip0.23in


There are many ways of defining the fractional Laplacian $(-\Delta)^s$, where $s$ is any positive, noninteger number.
Caffarelli and Silvestre in \cite{Caffarelli2007} gave a characterization of the fractional Laplacian when $0<s<1$ as the
Dirichlet-to-Neumann map for a function $u_e$ satisfying a higher order elliptic equation in the upper half space
with one extra spatial dimension. This idea  was  later generalized by Yang in \cite{Yang2100} when  the $s$ is being any positive, noninteger number.
See also Chang-Gonzales \cite{Chang2011} and Case-Chang \cite{Case2015} for general manifolds.

\vskip0.123in
To introduce the fractional operator $(-\Delta)^s$ for $2<s<3$, just like the case of $1<s<2$, via  the
Fourier transform, we can define
\be\nonumber
\widehat{(-\Delta)^s}u(\xi)=|\xi|^{2s} \widehat{u}(\xi)
\ee
or equivalently define this operator inductively by $(-\Delta)^s=(-\Delta)^{s-2}\circ(-\Delta)^2$.

\vskip0.1in
\begin{definition}
 We say a solution $u$ of \eqref{8LE} is stable outside a compact set if there exists $R_0>0$ such that
 \be\label{8stable}
 \int_{\R^n}\int_{\R^n}\frac{(\varphi(x)-\varphi(y))^2}{|x-y|^{n+2s}}dx dy-p\int_{\R^n}|u|^{p-1}\varphi^2dx\geq0
 \ee
 for any $\varphi\in C_c^{\infty}(\R^n\backslash \overline{B_{R_0}})$.
 \end{definition}
Set
\be\nonumber
p_s(n)=
\begin{cases}
\infty\;\;&\hbox{if}\;\; n\leq 2s,\\
\frac{n+2s}{n-2s}\;\;&\hbox{if}\;\; n>2s.
\end{cases}
\ee
The first main result  of the present paper is the following


\bt\label{8Liouvillec}
Suppose that $n>2s$ and $2<s<\delta<3$. Let $u\in C^{2\delta}(\R^n)\cap L^1(\R^n,(1+|z|)^{n+2s}dz)$ be a solution of
\eqref{8LE} which  is stable outside a compact set. Assume

\begin{itemize}
\item [(1)]   $1<p<p_s(n)$
or
\item [(2)]  $p_s(n)< p$  and
\be\label{8gamma}
p\frac{\Gamma(\frac{n}{2}-\frac{s}{p-1})\Gamma(s+\frac{s}{p-1})}{\Gamma(\frac{s}{p-1})\Gamma(\frac{n-2s}{2}-\frac{s}{p-1})}
>\frac{\Gamma(\frac{n+2s}{4})^2}{\Gamma(\frac{n-2s}{4})^2},
\ee
then the solution $u\equiv0$.
\item [(3)] $p=p_s(n)$, then $u$ has finite energy, i.e.,
\be\nonumber
\|u\|^2_{\dot{H}^s(\R^n)}:= \int_{\R^n}\int_{\R^n}\frac{(u(x)-u(y))^2}{|x-y|^{n+2s}}dx dy=\int_{\R^n}|u|^{p+1}<+\infty.
\ee
If in addition $u$ is stable, then $u\equiv0$.
\end{itemize}
\et


\br
In the Theorem  \ref{8Liouvillec} the condition \eqref{8gamma} is optimal. In fact,  the radial singular solution $u=|x|^{-\frac{2s}{p-1}}$ is stable if  \be\nonumber
p\frac{\Gamma(\frac{n}{2}-\frac{s}{p-1})\Gamma(s+\frac{s}{p-1})}{\Gamma(\frac{s}{p-1})\Gamma(\frac{n-2s}{2}-\frac{s}{p-1})}
\leq\frac{\Gamma(\frac{n+2s}{4})^2}{\Gamma(\frac{n-2s}{4})^2}.
\ee See \cite{LWZ}.
\er


\br The hypothesis $(2)$  of Theorem \ref{8Liouvillec} is  equivalent to
\be
p<p_c(n):=\begin{cases}
+\infty\;\;\;\;\;\;\;\;&\hbox{if}\;\;\;\;\;\;\;\; n\leq n_0(s),\\
\frac{n+2s-2-2a_{n,s}\sqrt{n}}{n-2s-2-2a_{n,s}\sqrt{n}}&\hbox{if}\;\;\;\;\;\;\;\;n> n_0(s),\\
\end{cases}
\ee
where $n_0(s)$ is the largest root  of $n-2s-2-2a_{n,s}\sqrt{n}=0$, see \cite{LWZ}. More details and further sharp results about $a_{n,s}$ and $n_0(s)$ see \cite{LWZ=00}.

\er


\br
In this remark, we further analyze the  hypothesis $(2)$  in Theorem \ref{8Liouvillec}. Recall that when $s=1$
the condition \eqref{8gamma} gives a upper bounded of $p$ (originated from Joseph and Lundgren \cite{Joseph1972}), it is
\be
p<p_c(n):=\begin{cases}
\;\;\;\;\;\;\;\;\;\infty\;\;\;\;\;\;\;\;&\hbox{if}\;\;\;\;\;\;\;\; n\leq10,\\
\frac{(n-2)^2-4n+8\sqrt{n-1}}{(n-2)(n-10)}\;\;\;\;\;\;\;\;&\hbox{if}\;\;\;\;\;\;\;\; n\geq11.
\end{cases}
\ee
For the case $s=2$, \eqref{8gamma} induce the upper bound of $p$ which is given by the following formula (cf. Gazzola and Grunau \cite{Gazzola2006}):
\be
p<p_c(n)=\begin{cases}
\;\;\;\;\;\;\;\;\;\;\infty\;\;\;\;\;\;\;\;&\hbox{if}\;\;\;\;\;\;\;\; n\leq12,\\
\frac{n+2-\sqrt{n^2+4-n\sqrt{n^2-8n+32}}}{n-6-\sqrt{n^2+4-n\sqrt{n^2-8n+32}}}\;\;\;\;\;\;\;\;&\hbox{if}\;\;\;\;\;\;\;\; n\geq13.
\end{cases}
\ee


\vskip0.1in

In the triharmonic case, the corresponding exponent given by see (\cite{LWZ2016=3}) is the following

\be\nonumber
p<p_{c}(n)=
\begin{cases}
\;\;\;\infty\;\;&\hbox{if}\;\; n\leq 14,\\
\frac{n+4-2D(n)}{n-8-2D(n)}\;\;&\hbox{if}\;\; n\geq15,
\end{cases}
\ee
where
\be\nonumber
D(n):=\frac{1}{6}\Big(9n^2+96-\frac{1536+1152n^2}{D_0(n)}-\frac{3}{2}D_0(n)\Big)^{1/2};
\ee
\be\nonumber
D_0(n):=-(D_1(n)+36\sqrt{D_2(n)})^{1/3};
\ee
\be\nonumber\aligned
D_1(n):=-94976+20736n+103104n^2-10368n^3+1296n^5-3024n^4-108n^6;
\endaligned\ee
\be\nonumber\aligned
D_2(n):&=6131712-16644096n^2+6915840n^4-690432n^6-3039232n\\
&\quad+4818944n^3-1936384n^5+251136n^7-30864n^8-4320n^9\\
&\quad+1800n^{10}-216n^{11}+9n^{12}.
\endaligned\ee

\er


\s{Preliminary}

Throughout this paper we denote $b:=5-2s$ and define the operator
\be\nonumber
\Delta_b w:=\Delta w+\frac{b}{y}w_y=y^{-b}\mathbf{div}(y^b\nabla w)
\ee
for a function $w\in W^{3,2}(\R^{n+1};y^bdxdy)$.  We firstly quote the following result.

\bt\label{ththYang} (See \cite{Yang2100} ) Assume $2<s<3$.  Let  $u_e\in W^{3,2}(\R^{n+1};y^bdxdy)$ satisfy  the equation
\be\label{8eqs}
\Delta_b^3u_e=0
\ee
on the upper half space for $(x,y)\in \R^n\times\R_+$ (where $y$ is the spacial direction) and the boundary conditions:
\be\label{8eqs1}\aligned
&u_e(x,0)=f(x),\\
&\lim_{y\rightarrow0}y^b\frac{\pa u_e}{\pa y}=0,\\
&\frac{\pa^2 u_e}{\pa y^2}\mid_{y=0}=\frac{1}{2s}\Delta_{x}u_e\mid_{y=0},\\
&\lim_{y\rightarrow0}C_{n,s}y^b\frac{\pa}{\pa y}\Delta_b^2 u_e=(-\Delta)^s f(x),
\endaligned
\ee
where $f(x)$ is some function defined on $H^s(\R^n)$.  Then we have
\be\label{Yangmore}
\int_{\R^n}|\xi|^{2s}|\hat{u}(\xi)|^2d\xi=C_{n,s}\int_{\R^{n+1}_+} y^b|\nabla\Delta_b u_e(x,y)|^2dxdy.
\ee
\et

\vskip0.23in

Applying the above theorem to solutions of \eqref{8LE},  we conclude that the extended function $u_e(x,y)$,
where $x=(x_1,...,x_n)\in\R^n$ and $y\in\R^+$, satisfies

\be\label{8LEE}\begin{cases}
\Delta_b^3 u_e=0\;\;\hbox{in}\;\;\R^{n+1}_+,\\
\lim_{y\rightarrow0}y^b\frac{\pa u_e}{\pa y}=0\;\;\hbox{on}\;\;\pa\R^{n+1}_+,\\
\frac{\pa^2 u_e}{\pa y^2}\mid_{y=0}=\frac{1}{2s}\Delta_{x}u_e\mid_{y=0}\;\;\hbox{on}\;\;\pa\R^{n+1}_+,\\
\lim_{y\rightarrow0}y^b\frac{\pa}{\pa y}\Delta_b^2 u_e=-C_{n,s}|u_e|^{p-1}u_e\;\;\hbox{in}\;\;\R^{n+1}_+.
\end{cases}\ee
Moreover,
\be\nonumber
\int_{\R^n}|\xi|^{2s}|\hat{u}(\xi)|^2d\xi=C_{n,s}\int_{\R^{n+1}_+} y^b|\nabla\Delta_b u_e(x,y)|^2dxdy
\ee
and $u(x)=u_e(x,0)$.

\vskip0.11in

\vskip0.1in


Define
\be\label{8exu}\aligned
&E(\la,x,u_e)=\int_{\R^{n+1}_+\cap\pa B_1}\frac{1}{2}\theta_1^b|\nabla\Delta_b u_e^\la|^2-\frac{C_{n,s}}{p+1}\int_{\pa \R^{n+1}_+\cap B_1}|u_e^\la|^{p+1}\\
&+\sum_{0\leq i,j\leq 4,i+j\leq5}C_{i,j}^1\int_{\R^{n+1}_+\cap\pa B_1}\theta_1^b\la^{i+j}\frac{d^i u^\la_e}{d\la^i}\frac{d^j u^\la_e}{d\la^j}\\
&+\sum_{0\leq t,s\leq 2,t+s\leq3}C_{t,s}^2\int_{\R^{n+1}_+\cap\pa B_1}\theta_1^b\la^{t+s}\nabla_{S^n}\frac{d^t u^\la_e}{d\la^t}\nabla_{S^n}\frac{d^s u^\la_e}{d\la^s}\\
&+\sum_{0\leq l,k\leq 1,l+k\leq1}C_{l,k}^3\int_{\R^{n+1}_+\cap\pa B_1}\theta_1^b\la^{l+k}\Delta_{S^n}\frac{d^l u^\la_e}{d\la^l}\Delta_{S^n}\frac{d^k u^\la_e}{d\la^k}\\
&+(\frac{s}{p-1}+1)\int_{\R^{n+1}_+\cap\pa B_1}\theta_1^b(\Delta_b u^\la_e)^2.
\endaligned\ee



\vskip0.1in

The following is the monotonicity formula which will paly an important role.
\bt\label{8monoid}
Let $u_e$ satisfy  the equation \eqref{8eqs} with the boundary conditions \eqref{8eqs1}, we have the following
\be\aligned
\frac{d E(\la,x,u_e)}{d\la}&=\int_{\R^{n+1}_+\cap\pa B_1}\theta_1^b\Big(
3\la^5(\frac{d^3u_e^\la}{d\la^3})^2+A_1\la^3(\frac{d^2u_e^\la}{d\la^2})^2+A_2\la(\frac{du^\la_e}{d\la})^2\Big)\\
&+\int_{\R^{n+1}_+\cap\pa B_1}\theta_1^b\Big(
2\la^3|\nabla_{S^n}\frac{d^2u_e^\la}{d\la^2}|^2+B_1\la|\nabla_{S^n}\frac{du_e^\la}{d\la}|\Big)\\
&+\int_{\R^{n+1}_+\cap\pa B_1}\theta_1^b\la(\Delta_{S^n}\frac{du_e^\la}{d\la})^2,
\endaligned
\ee
\et
where $\theta_1=\frac{y}{r}$ and
\be\nonumber\aligned
&A_1:=10\delta_1-2\delta_2-56+\alpha_0^2-2\al_0-2\beta_0-4,\\
&A_2:=-18\delta_1+6\delta_2-4\delta_3+2\delta_4+72-\al_0^2+\beta_0^2+2\al_0+2\beta_0,\\
&B_1:=8\al-4\beta-2\beta_0+4(n+b)-14,
\endaligned\ee
\be\nonumber\aligned
&\al:=n+b-2-\frac{4s}{p-1},\beta:=\frac{2s}{p-1}(3+\frac{2s}{p-1}-n-b),\quad\quad \\
&\al_0:=n+b-\frac{4s}{p-1},\;\;\;\beta_0:=\frac{2s}{p-1}(1+\frac{2s}{p-1}-n-b)\quad\quad
\endaligned\ee
and
\be\label{8delta}\aligned
\delta_1=&2(n+b)-\frac{8s}{p-1},\\
\delta_2=&(n+b)(n+b-2)-(n+b)\frac{12s}{p-1}+\frac{12s}{p-1}(1+\frac{2s}{p-1}),\\
\delta_3=&-\frac{8s}{p-1}(1+\frac{2s}{p-1})(2+\frac{2s}{p-1})+2(n+b)\frac{6s}{p-1}(1+\frac{2s}{p-1})\\
           &-(n+b)(n+b-2)(1+\frac{4s}{p-1}),\\
\delta_4=&(3+\frac{2s}{p-1})(2+\frac{2s}{p-1})(1+\frac{2s}{p-1})\frac{2s}{p-1}\\
&\quad-2(n+b)(1+\frac{2s}{p-1})(2+\frac{2s}{p-1})\frac{2s}{p-1}\\
          &\quad+(n+b)(n+b-2)(2+\frac{2s}{p-1})\frac{2s}{p-1}.\\
\endaligned\ee

We will give the proof of Theorem   \ref{8monoid} in the next section. Now we would like to
state a consequent result of  Theorem   \ref{8monoid}.
Recall that  $E(\la,x,u_e)$, defined in \eqref{8exu},  can be divided into two parts:  the integral  over  the ball $B_\la$
and the terms on the boundary $\pa B_\la$.  We note  that in our blow-down analysis, the coefficients (including positive or negative, big or small) of the boundary terms can be estimated  in a unified way, therefore we may change some coefficients of the boundary terms  in $E(\la,x,u_e)$. After such a change,  we denote the new functional  by $E^c(\la,x,u_e)$.

\vskip0.1in
Define
\be\label{728-pm}
p_m(n):=\begin{cases}
+\infty\;\;\;\;\;\;\;\;&\hbox{if}\;\;\;\;\;\;\;\; n<2s+6+\sqrt{73},\\
\frac{5n+10s-\sqrt{15(n-2s)^2+120(n-2s)+370}}{5n-10s-\sqrt{15(n-2s)^2+120(n-2s)+370}}&\hbox{if}\;\;\;\;\;\;\;\;n\geq2s+6+\sqrt{73}.\\
\end{cases}
\ee
 We have the following
\bt\label{8Monotonem} Assume that $\frac{n+2s}{n-2s}<p<p_m(n)$, then $E^c(\la,x,u_e)$ is a nondecreasing function of $\la>0$. Furthermore,
\be\nonumber
\frac{d E^c(\la,x,u_e)}{d\la}\geq C(n,s,p)\la^{2s\frac{p+1}{p-1}-6-n}
\int_{\R^{n+1}_+\cap\pa B_\la(x_0)}y^b\Big(\frac{2s}{p-1}u_e+\la\pa_r u_e\Big)^2,
\ee
where $C(n,s,p)$ is a constant independent of $\la$.
\et

By carefully comparing $\frac{n+2s}{n-2s}<p<p_m(n)$ with $p>\frac{n+2s}{n-2s}$ and \eqref{8gamma},  we get the following (see the last section of the current paper) monotonicity formula for our blow down analysis.

\bt\label{8Monotone} Assume that $p>\frac{n+2s}{n-2s}$ and \eqref{8gamma}, then $E^c(\la,x,u_e)$ is a nondecreasing function of $\la>0$. Furthermore,
\be\nonumber
\frac{d E^c(\la,x,u_e)}{d\la}\geq C(n,s,p)\la^{2s\frac{p+1}{p-1}-6-n}
\int_{\R^{n+1}_+\cap\pa B_\la(x_0)}y^b\Big(\frac{2s}{p-1}u_e+\la\pa_r u_e\Big)^2,
\ee
where $C(n,s,p)$ is a constant independent of $\la$.
\et

\section{Monotonicity formula and the proof of Theorem   \ref{8monoid}}

The derivation of the monotonicity for the \eqref{8LE} when $2<s<3$ is complicated in its process, we divide it into several subsections.
In subsection $3.1$, we  derive  $\frac{d}{d\la}\overline{E}(u_e,\la)$.
In subsection $3.2$, we calculate  $\frac{\pa^j}{\pa r^j}u_e^\la$  and  $\frac{\pa^i}{\pa \la^i}u_e^\la,\;\; i,j=1,2,3,4$.
In subsection $3.3$,  the operator   $\Delta_b^2$   and its representation will be given.
In subsection $3.4$, we  decompose $\frac{d}{d\la}\overline{E}(u_e^\la,1)$.
Finally, combine with the above four subsections, we can obtain the  monotonicity formula, hence get the proof of  Theorem \ref{8monoid}.

\vskip0.1in

Suppose that $x_0=0$ and denote by $B_\la$  the balls  centered at zero with radius $\la$. Set
\be\nonumber
\overline{E}(u_e,\la):=\la^{2s\frac{p+1}{p-1}-n}\Big(\int_{\R^{n+1}_+\cap B_\la}\frac{1}{2}y^b|\nabla\Delta_b u_e|^2
-\frac{C_{n,s}}{p+1}\int_{\pa\R^{n+1}_+\cap B_\la}|u_e|^{p+1}\Big).\\
\ee

\subsection{The derivation of $\frac{d}{d\la}\overline{E}(u_e,\la)$}
Define
\be\label{8uvw}\aligned
&v_e:=\Delta_b u_e, u_e^\la(X):=\la^{\frac{2s}{p-1}}u_e(\la X),\;\; w_e(X):=\Delta_b v_e\\
&v_e^\la(X):=\la^{\frac{2s}{p-1}+2}v_e(\la X),\;\;  w_e^\la(X):=\la^{\frac{2s}{p-1}+4}w_e(\la X),
\endaligned\ee
where $X=(x,y)\in \R^{n+1}_+$.
Therefore,
\be\label{8lamda}
\Delta_b u_e^\la(X)=v_e^\la(X), \Delta_b v_e^\la(X)=w_e^\la(X).
\ee
Hence
\be\nonumber\aligned
&\Delta_b w_e^\la=0,\\
&\lim_{y\rightarrow0}y^b\frac{\pa u_e}{\pa y}=0,\\
&\frac{\pa^2 u_e}{\pa y^2}\mid_{y=0}=\frac{1}{2s}\Delta_x u_e\mid_{y=0},\\
&\lim_{y\rightarrow0}C_{n,s}y^b\frac{\pa}{\pa y}w_e^\la=-C_{n,s}|u_e|^{p-1}u_e.
\endaligned
\ee
In addition, differentiating \eqref{8lamda} with respect to $\la$ we have
\be\nonumber
\Delta_b \frac{du_e^{\la}}{d\la}=\frac{dv_e^{\la}}{d\la},\;\;\; \Delta_b \frac{dv_e^{\la}}{d\la}=\frac{dw_e^{\la}}{d\la}.
\ee
Note that
\be\nonumber
\overline{E}(u_e,\la)=\overline{E}(u_e^\la,1)=\int_{\R^{n+1}_+\cap B_1}\frac{1}{2}y^b|\nabla v_e^\la|^2
-\frac{C_{n,s}}{p+1}\int_{\pa\R^{n+1}_+\cap B_1}|u_e^\la|^{p+1}.
\ee
Taking derivative of the energy $\overline{E}(u_e^\la,1)$ with respect to $\la$ and integrating by part we have:

\be\label{8overE}\aligned
&\frac{d\overline{E}(u_e^\la,1)}{d\la}=\int_{\R^{n+1}_+\cap B_1}y^b\nabla v_e^\la\nabla\frac{dv_e^\la}{d\la}
-C_{n,s}\int_{\pa\R^{n+1}_+\cap B_1}|u_e^\la|^{p-1}u_e^\la\frac{du_e^\la}{d\la}\\
=&\int_{\pa(\R^{n+1}_+\cap B_1)}y^b\frac{\pa v_e^\la}{\pa n}\frac{dv_e^\la}{d\la}-
\int_{\R^{n+1}_+\cap B_1}(y^b\Delta v_e^\la+by^{b-1}\frac{\pa v_e^\la}{\pa y})\frac{dv_e^\la}{d\la}\\
&-C_{n,s}\int_{\pa\R^{n+1}_+\cap B_1}|u_e^\la|^{p-1}u_e^\la\frac{du_e^\la}{d\la}\\
=&-\int_{\pa \R^{n+1}\cap B_1}y^b\frac{\pa v_e^\la}{\pa y}\frac{dv_e^\la}{d\la}
+\int_{\R^{n+1}_+\cap\pa B_1}y^b\frac{\pa v_e^\la}{\pa r}\\
&-\int_{\R^{n+1}_+\cap B_1}y^b\Delta v_e^\la\frac{dv_e^\la}{d\la}+by^{b-1}\frac{\pa v_e^\la}{\pa y}\frac{v_e^\la}{d\la}
-\int_{\pa\R^{n+1}_+\cap B_1}y^b\frac{w_e^\la}{\pa y}\frac{u_e^\la}{d\la}.
\endaligned\ee
Now note that from the definition of $v_e^\la$ and by differentiating  it  with respect to $\la$,  we get the following identity for $X\in\R^{n+1}_+$,
\be\nonumber
r\frac{\pa v_e^\la}{\pa r}=\la \pa_{\la} v_e^\la-(\frac{2s}{p-1}+2)v_e^\la.
\ee
Hence,
\be\nonumber
\aligned
\int_{\R^{n+1}_+\cap B_1}y^b\frac{\pa v_e^\la}{\pa r}&\frac{dv_e^\la}{d\la}
=\int_{\R^{n+1}_+\cap B_1}y^b\Big(\la\frac{dv_e^\la}{d\la}\frac{dv_e^\la}{d\la}-(\frac{2s}{p-1}+2)v_e^\la\frac{dv_e^\la}{d\la}\Big)\\
&=\la\int_{\R^{n+1}_+\cap \pa B_1}y^b(\frac{dv_e^\la}{d\la})^2-(\frac{s}{p-1}+1)\frac{d}{d\la}\int_{\R^{n+1}_+\cap \pa B_1}y^b(v_e^\la)^2.\\
\endaligned
\ee
Note that
\be\nonumber
\aligned
-\int_{\R^{n+1}_+\cap  B_1}y^b\Delta v_e^\la\frac{dv_e^\la}{d\la}=&
\int_{\pa\R^{n+1}_+\cap  B_1}y^b\frac{\pa v_e^\la}{\pa y}\frac{v_e^\la}{d\la}-\int_{\R^{n+1}_+\cap \pa B_1}y^b\frac{\pa v_e^\la}{\pa r}\frac{v_e^\la}{d\la}\\
&+\int_{\R^{n+1}_+\cap  B_1}\nabla v_e^\la\nabla(y^b\frac{dv_e^\la}{d\la}).
\endaligned
\ee
Integration by part we have
\be\nonumber\aligned
&\int_{\R^{n+1}_+\cap  B_1}y^b\nabla v_e^\la \nabla\frac{dv_e^\la}{d\la}\\
&=-\int_{\pa\R^{n+1}_+\cap  B_1}y^b\frac{\pa v_e^\la}{\pa y}\frac{dv_e^\la}{d\la}+\int_{\R^{n+1}_+\cap  \pa B_1}y^b\frac{\pa v_e^\la}{\pa r}\frac{dv_e^\la}{d\la}\\
&\quad -\int_{\R^{n+1}_+\cap  B_1}\nabla\cdot(y^b\nabla v^\la_e)\frac{dv_e^\la}{d\la}\\
&=-\int_{\pa\R^{n+1}_+\cap  B_1}y^b\frac{\pa v_e^\la}{\pa y}\frac{dv_e^\la}{d\la}+\int_{\R^{n+1}_+\cap  \pa B_1}y^b\frac{\pa v_e^\la}{\pa r}\frac{dv_e^\la}{d\la}\\
&\quad -\int_{\R^{n+1}_+\cap  B_1}y^b\Delta_b v_e^\la\Delta_b\frac{du_e^\la}{d\la}.\\
\endaligned
\ee
Now
\be\nonumber\aligned
-\int_{\R^{n+1}_+\cap  B_1}&y^b\Delta_b v_e^\la\Delta_b\frac{du_e^\la}{d\la}
=-\int_{\R^{n+1}_+\cap  B_1}y^b\Delta_b v_e^\la(\Delta\frac{du_e\la}{d\la}+\frac{b}{y}\frac{\pa}{\pa y}\frac{du_e^\la}{d\la})\\
=&-\int_{\pa(\R^{n+1}_+\cap  B_1)}y^b\Delta_b v_e^\la\frac{\pa}{\pa n}\frac{du_e^\la}{d\la}
+\int_{\R^{n+1}_+\cap  B_1}\nabla(y^b\Delta_b v_e^\la)\nabla\frac{du_e^\la}{d\la}\\
&-\int_{\R^{n+1}_+\cap  B_1}by^{b-1}\Delta_b v_e^\la\frac{\pa}{\pa y} \frac{du_e^\la}{d\la}\\
=&-\int_{\pa(\R^{n+1}_+\cap  B_1)}y^b\Delta_b v_e^\la \frac{\pa}{\pa n}\frac{d u_e^\la}{d\la}
+\int_{\R^{n+1}_+\cap  B_1}y^b\nabla\Delta_b v_e^\la\nabla \frac{du_e^\la}{d\la}\\
=&-\int_{\pa(\R^{n+1}_+\cap  B_1)}y^b\Delta_b v_e^\la \frac{\pa}{\pa n}\frac{d u_e^\la}{d\la}
+\int_{\pa(\R^{n+1}_+\cap  B_1)}y^b\frac{\pa \Delta_b v_e^\la}{\pa n}\frac{d u_e^\la}{d\la}\\
&-\int_{\R^{n+1}_+\cap  B_1}y^b\Delta_b^2 v_e^\la\frac{d u_e^\la}{d\la}\\
=&-\int_{\pa(\R^{n+1}_+\cap  B_1)}y^b\Delta_b v_e^\la \frac{\pa}{\pa n}\frac{d u_e^\la}{d\la}
+\int_{\pa(\R^{n+1}_+\cap  B_1)}y^b\frac{\pa \Delta_b v_e^\la}{\pa n}\frac{d u_e^\la}{d\la}.\\
\endaligned
\ee
Here we have used that $\Delta_b^2 v_e^\la=\Delta_b^3 u_e^\la=0$.
Therefore, combine with the above arguments we get that
\be\label{8overE1}\aligned
\int_{\R^{n+1}_+\cap  B_1}&y^b\nabla v_e^\la \nabla\frac{dv_e^\la}{d\la}=-\int_{\pa\R^{n+1}_+\cap  B_1}y^b\frac{\pa v_e^\la}{\pa y}\frac{dv_e^\la}{d\la}
+\int_{\R^{n+1}_+\cap \pa B_1}y^b \frac{\pa v_e^\la}{\pa r}\frac{dv_e^\la}{d\la}\\
&+\int_{\pa\R^{n+1}_+\cap  B_1}y^b\Delta_b v_e^\la\frac{\pa}{\pa y}\frac{du_e^\la}{d\la}-\int_{\R^{n+1}_
+\cap \pa B_1}y^b \Delta v_e^\la \frac{\pa}{\pa r}\frac{du_e^\la}{d\la}\\
&-\int_{\pa\R^{n+1}_+\cap  B_1}y^b\frac{\pa}{\pa y}\Delta_b v_e^\la\frac{du_e^\la}{d\la}+\int_{\R^{n+1}_+\cap \pa B_1}y^b\frac{\pa}{\pa r}\Delta_b v_e^\la\frac{du_e^\la}{d\la}\\
=&\int_{\R^{n+1}_+\cap \pa B_1}y^b \frac{\pa v_e^\la}{\pa r}\frac{dv_e^\la}{d\la}-\int_{\R^{n+1}_
+\cap \pa B_1}y^b \Delta v_e^\la \frac{\pa}{\pa r}\frac{du_e^\la}{d\la}\\
&-C(n,s)\int_{\pa\R^{n+1}_+\cap  B_1}|u_e^\la|^{p-1}\frac{du_e^\la}{d\la}+\int_{\R^{n+1}_+\cap \pa B_1}y^b\frac{\pa}{\pa r}\Delta_b v_e^\la\frac{du_e^\la}{d\la}.\\
\endaligned
\ee
Here,  we have used that $\frac{\pa \Delta u_e^\la(x,0)}{\pa y}=0$, $\frac{\pa}{\pa y}\frac{du_e^\la}{d\la}=0$ on $\pa \R^{n+1}_+$
 and $\lim_{y\rightarrow0}y^b\frac{\pa}{\pa y}\Delta_b v_e^\la=-C_{n,s}|u_e^\la|^{p-1}u_e^\la$.
By \eqref{8overE} and \eqref{8overE1} we obtain that
\be\label{8overee}\aligned
\frac{d}{d\la}\overline{E}(u_e^\la,1)=&\int_{\R^{n+1}_+\cap \pa B_1}y^b\frac{\pa v_e^\la}{\pa r}\frac{dv_e^\la}{d\la}
+\int_{\R^{n+1}_+\cap \pa B_1}y^b\frac{\pa w_e^\la}{\pa r}\frac{du_e^\la}{d\la}\\
&-\int_{\R^{n+1}_+\cap \pa B_1}y^bw_e^\la\frac{\pa}{\pa r}\frac{du_e^\la}{d\la}.
\endaligned\ee
Recall  \eqref{8uvw} and
differentiate it with respect to $\la$, we have
\be\nonumber\aligned
\frac{du_e^\la(X)}{d\la}&=\frac{1}{\la}\big(\frac{2s}{p-1}u_e^\la(X)+r\pa_r u_e^\la(X)\big),\\
\frac{dv_e^\la(X)}{d\la}&=\frac{1}{\la}\big((\frac{2s}{p-1}+2)v_e^\la(X)+r\pa_r v_e^\la(X)\big),\\
\frac{dw_e^\la(X)}{d\la}&=\frac{1}{\la}\big((\frac{2s}{p-1}+4)w_e^\la(X)+r\pa_r w_e^\la(X)\big).\\
\endaligned
\ee
Differentiate the above equations with respect to $\la$ again we get
\be\nonumber
\la\frac{d^2u_e^\la(X)}{d\la^2}+\frac{du_e^\la(x)}{d\la}=\frac{2s}{p-1}\frac{du_e^\la(X)}{d\la}+r\pa_r\frac{du_e^\la}{d\la}.
\ee
Hence, for $X\in\R^{n+1}_+\cap B_1$, we have
\be\nonumber\aligned
\pa_r(u_e^\la(X))&=\la\frac{du_e^\la}{d\la}-\frac{2s}{p-1}u_e,\\
\pa_r(\frac{du_e^\la(X)}{d\la})&=\la\frac{d^2u_e^\la(X)}{d\la^2}+(1-\frac{2s}{p-1})\frac{du_e^\la}{d\la},\\
\pa_r(v_e^\la(X))&=\la\frac{dv_e^\la}{d\la}-(\frac{2s}{p-1}+2)v_e^\la,\\
\pa_r(w_e^\la(X))&=\la\frac{d w_e^\la}{d\la}-(\frac{2s}{p-1}+4)w_e^\la.\\
\endaligned
\ee
Plugging these equations into \eqref{8overee}, we get that
\be\label{8overlineE}\aligned
\frac{d}{d\la}\overline{E}(u_e^\la,1)=&\int_{\R^{n+1}_+\cap \pa B_1}y^b\big(\la\frac{dv_e^\la}{d\la}\frac{dv_e^\la}{d\la}-(\frac{2s}{p-1}+2)v_e^\la\frac{dv_e^\la}{d\la}\big)\\
&+y^b\big(\la\frac{d w_e^\la}{d\la}\frac{u_e^\la}{d\la}-(\frac{2s}{p-1}+4)w_e^\la\frac{du_e^\la}{d\la}\big)\\
&-y^b\big(\la w_e^\la\frac{d^2 u_e^\la}{d\la^2}+(1-\frac{2s}{p-1}w_e^\la\frac{du_e^\la}{d\la}\big)\\
 =&\underbrace{\int_{\R^{n+1}_+\cap \pa B_1}y^b\big[\la\frac{dv_e^\la}{d\la}\frac{dv_e^\la}{d\la}-(\frac{2s}{p-1}+2)v_e^\la\frac{dv_e^\la}{d\la}\big]}\\
&+\underbrace{y^b\big[\la \frac{d w_e^\la}{d\la}\frac{du_e^\la}{d\la}-\la w_e^\la\frac{d^2u_e^\la}{d\la^2}\big]-5y^bw_e^\la\frac{du_e^\la}{d\la}}\\
:=&\overline{E}_{d1}(u_e^\la,1)+\overline{E}_{d2}(u_e^\la,2).
\endaligned
\ee

\subsection{The calculations  of $\frac{\pa^j}{\pa r^j}u_e^\la$  and  $\frac{\pa^i}{\pa \la^i}u_e^\la,\; i,j=1,2,3,4$}

Note
\be\label{8ue0}
\la\frac{du_e^\la}{d\la}=\frac{2s}{p-1}u_e^\la+r\frac{\pa}{\pa r}u_e^\la.
\ee
Differentiating \eqref{8ue0} once, twice and  thrice with respect to $\la$ respectively,  we have
\be\label{8ue1}
\la\frac{d^2u_e^\la}{d\la^2}+\frac{du_e^\la}{d\la}=\frac{2s}{p-1}\frac{du_e^\la}{d\la}+r\frac{\pa}{\pa r}\frac{du_e^\la}{d\la},
\ee
\be\label{8ue2}
\la\frac{d^3u_e^\la}{d\la^3}+2\frac{d^2u_e^\la}{d\la^2}=\frac{2s}{p-1}\frac{d^2u_e^\la}{d\la^2}+r\frac{\pa}{\pa r}\frac{d^2u_e^\la}{d\la^2},
\ee
\be\label{8ue3}
\la\frac{d^4u_e^\la}{d\la^4}+3\frac{d^3u_e^\la}{d\la^3}=\frac{2s}{p-1}\frac{d^3u_e^\la}{d\la^3}+r\frac{\pa}{\pa r}\frac{d^3u_e^\la}{d\la^3}.
\ee
Similarly, differentiating \eqref{8ue0} once, twice and  thrice with respect to $r$ respectively we have
\be\label{8ue4}
\la\frac{\pa}{\pa r}\frac{du_e^\la}{d\la}=(\frac{2s}{p-1}+1)\frac{\pa}{\pa r}u_e^\la+r\frac{\pa^2}{\pa r^2}u_e^\la,
\ee
\be\label{8ue5}
\la\frac{\pa^2}{\pa r^2}\frac{du_e^\la}{d\la}=(\frac{2s}{p-1}+2)\frac{\pa^2}{\pa r^2}u_e^\la+r\frac{\pa^3}{\pa r^3}u_e^\la,
\ee
\be\label{8ue6}
\la\frac{\pa^3}{\pa r^3}\frac{du_e^\la}{d\la}=(\frac{2s}{p-1}+3)\frac{\pa^3}{\pa r^3}u_e^\la+r\frac{\pa^4}{\pa r^4}u_e^\la.
\ee
From \eqref{8ue0}, on $\R^{n+1}_+\cap \pa B_1$, we have
\be\nonumber
\frac{\pa u_e^\la}{\pa r}=\la \frac{du_e^\la}{d\la}-\frac{2s}{p-1}u_e^\la.
\ee
Next from \eqref{8ue1}, on $\R^{n+1}_+\cap \pa B_1$,  we derive that
\be\nonumber
\frac{\pa}{\pa r}\frac{d u_e^\la}{d\la}=\la\frac{d^2u_e^\la}{d\la^2}+(1-\frac{2s}{p-1})\frac{du_e^\la}{d\la}.
\ee
From \eqref{8ue4}, combine with the two equations above, on $\R^{n+1}_+\cap \pa B_1$, we get
\be\label{8r31}\aligned
\frac{\pa^2}{\pa r^2} u_e^\la&=\la\frac{\pa}{\pa r}\frac{du_e^\la}{d\la}-(1+\frac{2s}{p-1})\frac{\pa}{\pa r} u_e^\la\\
&=\la^2\frac{d^2 u_e^\la}{d\la^2}-\la\frac{4s}{p-1}\frac{du_e^\la}{d\la}+(1+\frac{2s}{p-1})\frac{2s}{p-1}u_e^\la.
\endaligned\ee
Differentiating \eqref{8ue1} with respect to $r$, and combine with \eqref{8ue1} and \eqref{8ue2}, we get that
\be\label{8r32}\aligned
\frac{\pa^2}{\pa r^2}\frac{du_e^\la}{d\la}&=\la\frac{\pa}{\pa r}\frac{d^2 u_e^\la}{d\la^2}-\frac{2s}{p-1}\frac{\pa}{\pa r}\frac{du_e^\la}{d\la}\\
&=\la^2\frac{d^3 u_e^\la}{d\la^3}+(2-\frac{4s}{p-1})\la\frac{d^2 u_e^\la}{d\la^2}-(1-\frac{2s}{p-1})\frac{2s}{p-1}\frac{du_e^\la}{d\la}.
\endaligned\ee
From \eqref{8ue5}, on $\R^{n+1}_+\cap \pa B_1$, combine with \eqref{8r31} and \eqref{8r32}, we have
\be\label{8ue7}\aligned
\frac{\pa^3}{\pa r^3}u_e^\la=&\la\frac{\pa^2}{\pa r^2}\frac{du_e^\la}{d\la}-(2+\frac{2s}{p-1})\frac{\pa^2}{\pa r^2}u_e^\la\\
=&\la^3\frac{d^3 u_e^\la}{d\la^3}-\la^2\frac{6s}{p-1}\frac{d^2 u_e^\la}{d\la^2}+\la(\frac{6s}{p-1}+\frac{12s^2}{(p-1)^2})\frac{du_e^\la}{d\la}\\
&-(2+\frac{2s}{p-1})(1+\frac{2s}{p-1})\frac{2s}{p-1}u_e^\la.
\endaligned\ee
Now differentiating \eqref{8ue1} once with respect to $r$, we get
\be\nonumber
\la\frac{\pa^2}{\pa r^2}\frac{d^2 u_e^\la}{d\la^2}=(\frac{2s}{p-1}+1)\frac{\pa^2}{\pa r^2}\frac{du_e^\la}{d\la}+r\frac{\pa^3}{\pa r^3}\frac{du_e^\la}{d\la},
\ee
then on $\R^{n+1}_+\cap \pa B_1$, we have
\be\label{8ue8}
\frac{\pa^3}{\pa r^3}\frac{du_e^\la}{d\la}=\la\frac{\pa^2}{\pa r^2}\frac{d^2 u_e^\la}{d\la^2}-(\frac{2s}{p-1}+1)\frac{\pa^2}{\pa r^2}\frac{du_e^\la}{d\la}.
\ee
Now differentiating \eqref{8ue2} twice with respect to $r$, we get
\be\nonumber
\la\frac{\pa}{\pa r}\frac{d^3 u_e^\la}{d\la^3}=(\frac{2s}{p-1}-1)\frac{\pa}{\pa r}\frac{d^2 u_e^\la}{d\la^2}+r\frac{\pa^2}{\pa r^2}\frac{d^2 u_e^\la}{d\la^2},
\ee
hence on $\R^{n+1}_+\cap \pa B_1$, combine with \eqref{8ue2} and \eqref{8ue3} there holds
\be\label{8ue9}
\aligned
\frac{\pa^2}{\pa r^2}&\frac{d^2 u_e^\la}{d\la^2}=\la\frac{\pa}{\pa r}\frac{d^3 u_e^\la}{d\la^3}+(1-\frac{2s}{p-1})\frac{\pa}{\pa r}\frac{d^2 u_e^\la}{d\la^2}\\
=&\la^2\frac{d^4 u_e^\la}{d\la^4}+\la(4-\frac{4s}{p-1})\frac{d^3 u_e^\la}{d\la^3}+(1-\frac{2s}{p-1})(2-\frac{2s}{p-1})\frac{d^2 u_e^\la}{d\la^2}.
\endaligned\ee
Now differentiating \eqref{8ue1} with respect to $r$, we have
\be\nonumber
\la\frac{\pa}{\pa r}\frac{d^2 u_e^\la}{d\la^2}=\frac{2s}{p-1}\frac{\pa}{\pa r}\frac{du_e^\la}{d\la}+r\frac{\pa^2}{\pa r^2}\frac{du_e^\la}{d\la}.
\ee
This combine with \eqref{8ue1} and \eqref{8ue2}, on $\R^{n+1}_+\cap \pa B_1$, we have
\be\label{8ue10}\aligned
\frac{\pa^2}{\pa r^2}\frac{du_e^\la}{d\la}&=\la\frac{\pa}{\pa r}\frac{d^2 u_e^\la}{d\la^2}-\frac{2s}{p-1}\frac{\pa}{\pa r}\frac{d u_e^\la}{d\la}\\
&=\la^2\frac{d^3 u_e^\la}{d\la^3}+\la(2-\frac{4s}{p-1})\frac{d^2 u_e^\la}{d\la^2}-\frac{2s}{p-1}(1-\frac{2s}{p-1})\frac{du_e^\la}{d\la}.
\endaligned\ee
Now from \eqref{8ue8}, combine with \eqref{8ue9} and \eqref{8ue10}, we get
\be\label{8ue11}\aligned
\frac{\pa^3}{\pa r^3}\frac{du_e^\la}{d\la}=&\la^3\frac{d^4 u_e^\la}{d\la^4}+\la^2(3-\frac{6s}{p-1})\frac{d^3 u_e^\la}{d\la^3}-\la(1-\frac{2s}{p-1})\frac{6s}{p-1}\frac{d^2 u_e^\la}{d\la^2}\\
&+(1-\frac{2s}{p-1})(1+\frac{2s}{p-1})\frac{2s}{p-1}\frac{du_e^\la}{d\la}.
\endaligned\ee
From \eqref{8ue6}, on $\R^{n+1}_+\cap \pa B_1$, combine with \eqref{8ue11} then
\be\nonumber\aligned
\frac{\pa^4}{\pa r^4}u_e^\la=&\la \frac{\pa^3}{\pa r^3}\frac{du_e^\la}{d\la}-(3+\frac{2s}{p-1})\frac{\pa^3}{\pa r^3}u_e^\la\\
=&\la^4\frac{d^4 u_e^\la}{d\la^4}-\la^3\frac{8s}{p-1}\frac{d^3 u_e^\la}{d\la^3}+\la^2(2+\frac{4s}{p-1})\frac{6s}{p-1}\frac{d^2 u_e^\la}{d\la^2}\\
&-\la (1+\frac{2s}{p-1})(1+\frac{s}{p-1})\frac{16s}{p-1}\frac{du_e^\la}{d\la}\\
&+(3+\frac{2s}{p-1})(2+\frac{2s}{p-1})(1+\frac{2s}{p-1})\frac{2s}{p-1}u_e^\la.\\
\endaligned\ee
In summary, we have that
\be\nonumber\aligned
\frac{\pa^3}{\pa r^3}u_e^\la
=&\la^3\frac{d^3 u_e^\la}{d\la^3}-\la^2\frac{6s}{p-1}\frac{d^2 u_e^\la}{d\la^2}+\la(\frac{6s}{p-1}+\frac{12s^2}{(p-1)^2})\frac{du_e^\la}{d\la}\\
&-(2+\frac{2s}{p-1})(1+\frac{2s}{p-1})\frac{2s}{p-1}u_e^\la
\endaligned\ee
and
\be\nonumber
\frac{\pa ^2}{\pa r^2}u_e^\la=\la^2\frac{d^2 u_e^\la}{d\la^2}-\la\frac{4s}{p-1}\frac{du_e^\la}{d\la}+(1+\frac{2s}{p-1})\frac{2s}{p-1}u_e^\la
\ee
\be\nonumber
\frac{\pa u_e^\la}{\pa r}=\la \frac{du_e^\la}{d\la}-\frac{2s}{p-1}u_e^\la.
\ee
\subsection{On the operator $\Delta_b^2$  and its representation }

Note that
\be\nonumber\aligned
\Delta_b u=&y^{-b}\nabla\cdot(y^b\nabla u)
=&u_{rr}+\frac{n+b}{r}u_r+\frac{1}{r^2}\theta_1^{-b}\mathbf{ div} _{S^n}(\theta_1^b\nabla_{S^n}u),
\endaligned\ee
where $\theta_1=\frac{y}{r},r=\sqrt{|x|^2+y^2}$. Set $v=\Delta_b u$  and $\Delta_b^2 u:=w$.
Then
\be\nonumber\aligned
w=&\Delta_b v=v_{rr}+\frac{n+b}{r}v_r+\frac{1}{r^2}\theta_1^{-b} \mathbf{ div} _{S^n}(\theta_1^b\nabla_{S^n}v)\\
=&\pa_{rrrr}u+\frac{2(n+b)}{r}\pa_{rrr}u+\frac{(n+b)(n+b-2)}{r^2}\pa_{rr}u-\frac{(n+b)(n+b-2)}{r^3}\pa_r u\\
&+r^{-4}\theta_1^{-b}{\bf div}_{S^n}(\theta_1^b\nabla(\theta_1^{-b}\mathbf{ div}_{S^n}(\theta_1^b\nabla_{S^n}u))\\
&+2r^{-2}\theta_1^{-b}{\bf div}_{S^n}(\theta_1^b\nabla_{S^n}(u_{rr}+\frac{n+b-2}{r}u_r))\\
&-2(n+b-3)r^{-4}\theta_1^{-b}\mathbf{ div}_{S^n}(\theta_1^b\nabla_{S^n}u).\\
\endaligned\ee
On $\R^{n+1}_+\cap \pa B_1$, we have
\be\nonumber\aligned
w=&\underbrace{\pa_{rrrr}u+2(n+b)\pa_{rrr}u+(n+b)(n+b-2)\pa_{rr}u-(n+b)(n+b-2)\pa_r u}\\
&\underbrace{+\theta_1^{-b}\mathbf{ div}_{S^n}(\theta_1^b\nabla(\theta_1^{-b}\mathbf{div}_{S^n}(\theta_1^b\nabla_{S^n}u))}\\
&\underbrace{+2\theta_1^{-b}\mathbf{div}_{S^n}(\theta_1^b\nabla_{S^n}(u_{rr}+\frac{n+b-2}{r}u_r))}\\
&\underbrace{-2(n+b-3)\theta_1^{-b}\mathbf{div}_{S^n}(\theta_1^b\nabla_{S^n}u)}\\
 :=&I(u)+J(u)+K(u)+L(u).\\
\endaligned\ee
By these notations, we can rewrite the term $\overline{E}_{d2}(u_e^\la,1)$ appear in \eqref{8overlineE} as following
\be\label{8ijk}\aligned
&\overline{E}_{d2}(u_e^\la,1)\\
&=\int_{\R^{n+1}_+\cap \pa B_1}\theta_1^b(\la\frac{d w_e^\la}{d\la}\frac{d u_e^\la}{d\la}-\la w_e^\la\frac{d^2 u_e^\la}{d\la^2})-5\theta_1^b w_e^\la\frac{d u_e^\la}{d\la}\\
&=\underbrace{\int_{\R^{n+1}_+\cap\pa B_1}\la\theta_1^b\frac{d}{d\la}I(u_e^\la)\frac{d u_e^\la}{d\la}-\la\theta_1^b I(u_e^\la)\frac{d^2 u_e^\la}{d\la^2}-5\theta_1^bI(u_e^\la)\frac{du_e^\la}{d\la}}\\
&\quad \underbrace{+\int_{\R^{n+1}_+\cap\pa B_1}\la\theta_1^b\frac{d}{d\la}J(u_e^\la)\frac{d u_e^\la}{d\la}-\la\theta_1^b J(u_e^\la)\frac{d^2 u_e^\la}{d\la^2}-5\theta_1^bJ(u_e^\la)\frac{du_e^\la}{d\la}}\\
&\quad \underbrace{+\int_{\R^{n+1}_+\cap\pa B_1}\la\theta_1^b\frac{d}{d\la}K(u_e^\la)\frac{d u_e^\la}{d\la}-\la\theta_1^b K(u_e^\la)\frac{d^2 u_e^\la}{d\la^2}-5\theta_1^bK(u_e^\la)\frac{du_e^\la}{d\la}}\\
&\quad \underbrace{+\int_{\R^{n+1}_+\cap \pa B_1}\la\theta_1^b\frac{d}{d\la}L(u_e^\la)\frac{du_e^\la}{d\la}-\la\theta_1^bL(u_e^\la)\frac{d^2u_e^\la}{d\la^2}-5\theta_1^bL(u_e^\la)\frac{du_e^\la}{d\la}},
\endaligned\ee
we define as
\be\nonumber\aligned
\overline{E}_{d2}(u_e^\la,1):&=\mathcal{I}+\mathcal{J}+\mathcal{K}+\mathcal{L}\\
:&=I_1+I_2+I_3+J_1+J_2+J_3+K_1+K_2+K_3+L_1+L_2+L_3.
\endaligned\ee
where $I_1,I_2,I_3,J_1,J_2,J_3,K_1,K_2,K_3,L_1,L_2,L_3 $ are corresponding successively to the $12$ terms in \eqref{8ijk}.
By the conclusions of subsection $2.2$, we have
\be\label{8iuiu}\aligned
I(u_e^\la)&=\pa_{rrrr} u_e^\la+2(n+b)\pa_{rrr} u_e^\la\\
&\quad +(n+b)(n+b-2)\pa_{rr} u_e^\la-(n+b)(n+b-2)\pa_r u_e^\la\\
&=\la^4\frac{d^4 u_e^\la}{d\la^4}+\la^3\big(2(n+b)-\frac{8s}{p-1})\big)\frac{d^3 u_e^\la}{d\la^3}\\
&\quad+\la^2\big[\frac{12s}{p-1}(1+\frac{2s}{p-1})-(n+b)\frac{12s}{p-1}+(n+b)(n+b-2)\big]\frac{d^2 u_e^\la}{d\la^2}\\
&\quad+\la\big[-\frac{8s}{p-1}(1+\frac{2s}{p-1})(2+\frac{2s}{p-1})+2(n+b)\frac{6s}{p-1}(1+\frac{2s}{p-1})\\
&\quad+(n+b)(n+b-2)(-\frac{4s}{p-1}-1)\big]\frac{d u_e^\la}{d\la}\\
&\quad+\big[(1+\frac{2s}{p-1})(2+\frac{2s}{p-1})(3+\frac{2s}{p-1})\frac{2s}{p-1}\\
&\quad-(n+b)(1+\frac{2s}{p-1})(2+\frac{2s}{p-1})\frac{4s}{p-1}\\
&\quad+(n+b)(n+b-2)(\frac{2s}{p-1}+2)\frac{2s}{p-1}\big]u_e^\la.\\
\endaligned\ee
For convenience, we denote that
\be\label{8iue}
I(u_e^\la)=\la^4\frac{d^4 u_e^\la}{d\la^4}+\la^3 \delta_1\frac{d^3 u_e^\la}{d\la^3}+\la^2\delta_2\frac{d^2 u_e^\la}{d\la^2}+\la\delta_3\frac{du_e^\la}{d\la}+\delta_4 u_e^\la,
\ee
where $\delta_i$ are the corresponding coefficients of $\la^i\frac{d^i u_e^\la}{d\la^i}$ appeared in \eqref{8iuiu} for $i=1,2,3,4$.
Now taking the derivative of \eqref{8iue} with respect to $\la$, we get
\be\label{8diue}\aligned
\frac{d}{d\la}I(u_e^\la)=&\la^4\frac{d^5 u_e^\la}{d\la^5}+\la^3 (\delta_1+4)\frac{d^4 u_e^\la}{d\la^4}+\la^2(3\delta_1+\delta_2)\frac{d^3 u_e^\la}{d\la^3}\\
&+\la(2\delta_2+\delta_3)\frac{d^2u_e^\la}{d\la^2}+(\delta_3+\delta_4)\frac{du_e^\la}{d\la}
\endaligned\ee
and
\be\label{8urr}\aligned
\pa_{rr}& u_e^\la+(n+b-2)\pa_r u_e^\la \\ =&\la^2\frac{d^2u_e^\la}{d\la^2}+\la(n+b-2-\frac{4s}{p-1})\frac{du_e^\la}{d\la}+\frac{2s}{p-1}(3+\frac{2s}{p-1}-n-b)u_e^\la\\
:=&\la^2\frac{d^2 u_e^\la}{d\la^2}+\la\alpha\frac{du_e^\la}{d\la}+\beta u_e^\la.
\endaligned\ee
Hence,
\be\label{8urrla}\aligned
\frac{d}{d\la}[\pa_{rr}& u_e^\la+(n+b-2)\pa_r u_e^\la]=\la^2\frac{d^3 u_e^\la}{d\la^3}+\la(\alpha+2)\frac{d^2 u_e^\la}{d\la^2}+(\alpha+\beta)\frac{du_e^\la}{d\la},
\endaligned\ee
here $\alpha=n+b-2-\frac{4s}{p-1}$ and $\beta=\frac{2s}{p-1}(3+\frac{2s}{p-1}-n-b)$.

\subsection{The computations  of $I_1,I_2,I_3$ and $\mathcal{I}$}

\be\label{8i1}\aligned
I_1:=&\int_{\R^{n+1}_+\cap\pa B_1}\la \theta_1^b\frac{d}{d\la}I(u_e^\la)\frac{du_e^\la}{d\la}\\
=&\int_{\R^{n+1}_+\cap\pa B_1}\theta_1^b\big(\la^5\frac{d^5 u_e^\la}{d\la^5}+\la^4(4+\delta_1)\frac{d^4 u_e^\la}{d\la^4}+\la^3(3\delta_1+\delta_2)\frac{d^3 u_e^\la}{d\la^3}\\
&+\la^2(2\delta_2+\delta_3)\frac{d^2 u_e^\la}{d\la^2}+\la(\delta_3+\delta_4)\frac{du_e^\la}{d\la}\big)\frac{du_e^\la}{d\la}\\
=&\frac{d}{d\la}\int_{\R^{n+1}_+\cap\pa B_1}\theta_1^b\big[\la^5\frac{d^4 u_e^\la}{d\la^4}\frac{d u_e^\la}{d\la}-\la^5\frac{d^3 u_e^\la}{d\la^3}\frac{d^2 u_e^\la}{d\la^2}+(\delta_1-1)\la^4\frac{d^3 u_e^\la}{d\la^3}\frac{d u_e^\la}{d\la}\\
&+(4-\delta_1+\delta_2)\la^3\frac{d^2 u_e^\la}{d\la^2}\frac{d u_e^\la}{d\la}+\frac{3\delta_1-\delta_2+\delta_3-12}{2}\la^2(\frac{d u_e^\la}{d\la})^2\big]\\
&+\int_{\R^{n+1}_+\cap\pa B_1}\theta_1^b \big[(12-3\delta_1+\delta_2+\delta_4)\la(\frac{d u_e^\la}{d\la})^2\\
&+(\delta_1-4-\delta_2)\la^3(\frac{d^2 u_e^\la}{d\la^2})^2+\la^5(\frac{d^3 u_e^\la}{d\la^3})^2\big]\\
&+\int_{\R^{n+1}_+\cap\pa B_1}\theta_1^b(6-\delta_1)\la^4\frac{d^3 u_e^\la}{d\la^3}\frac{d^2 u_e^\la}{d\la^2},
\endaligned
\ee
where $\delta_i (i=1,2,3,4) $ are defined in \eqref{8iuiu} and \eqref{8iue}. Denote   $f=u_e^\la,f':=\frac{du_e^\la}{d\la}$,   we have  used  the following differential identities:
\be\nonumber\aligned
\la^5 f'''''f'=&\big[\la^5f''''f'-\la^5f'''f''-5\la^4f'''f'+20\la^3f''f'-30\la^2f'f'\big]'\\
&+60\la(f')^2-20\la^3(f'')^2+\la^5(f''')^2+10\la^4f'''f'',
\endaligned\ee
\be\nonumber\aligned
\la^4f''''f'=\big[\la^4f'''f'-4\la^3f''f'+6\la^2f'f'\big]'-12\la(f')^2+4\la^3(f'')^2-\la^4f'''f'',
\endaligned\ee
\be\nonumber\aligned
\la^3f'''f'=\big[\la^3f''f'-\frac{3\la^2}{2}f'f'\big]'+3\la(f')^2-\la^3(f'')^2,
\endaligned\ee
and
\be\nonumber\aligned
\la^2 f''f'=\big[\frac{\la^2}{2}f'f'\big]'-\la(f')^2.
\endaligned\ee

\be\label{8i2}\aligned
&I_2:=-\la\int_{\R^{n+1}_+\cap\pa B_1}\theta_1^b I(u_e^\la)\frac{d^2 u_e^\la}{d\la^2}\\
&=-\la\int_{\R^{n+1}_+\cap\pa B_1}\theta_1^b \big(\la^4\frac{d^4 u_e^\la}{d\la^4}+\la^3 \delta_1\frac{d^3 u_e^\la}{d\la^3}+\la^2\delta_2\frac{d^2 u_e^\la}{d\la^2}+\la\delta_3\frac{du_e^\la}{d\la}+\delta_4 u_e^\la\big)\frac{d^2 u_e^\la}{d\la^2}\\
&=\frac{d}{d\la}\int_{\R^{n+1}_+\cap\pa B_1}\theta_1^b\big[-\la^5\frac{d^3 u_e^\la}{d\la^3}\frac{d^2 u_e^\la}{d\la^2}-\delta_4\la\frac{d u_e^\la}{d\la}u_e^\la\big]\\
&\quad+\int_{\R^{n+1}_+\cap\pa B_1}\theta_1^b\big[\la^5(\frac{d^3 u_e^\la}{d\la^3})^2-\delta_2\la^3(\frac{d^2 u_e^\la}{d\la^2})^2+\delta_4\la(\frac{d u_e^\la}{d\la})^2\big]\\
&\quad +\int_{\R^{n+1}_+\cap\pa B_1}\theta_1^b\big[(5-\delta_1)\la^4\frac{d^3 u_e^\la}{d\la^3}\frac{d^2 u_e^\la}{d\la^2}-\delta_3\la^2\frac{d^2 u_e^\la}{d\la^2}\frac{d u_e^\la}{d\la}+\delta_4\frac{d u_e^\la}{d\la}u_e^\la\big].
\endaligned\ee
Here we have used that
\be\nonumber\aligned
-\la^5f''''f''=\big[-\la^5f'''f''\big]'+5\la^4f'''f''+\la^5(f''')^2
\endaligned\ee
and
\be\nonumber\aligned
-\la f''f=[-\la f'f]'+f'f+\la(f')^2.
\endaligned\ee

\be\label{8iu3}\aligned
I_3:=&-5\int_{\R^{n+1}_+\cap\pa B_1}\theta_1^b I(u_e^\la)\frac{d u_e^\la}{d\la}\\
=&-5\int_{\R^{n+1}_+\cap\pa B_1}\theta_1^b\big[\la^4\frac{d^4 u_e^\la}{d\la^4}+\la^3 \delta_1\frac{d^3 u_e^\la}{d\la^3}+\la^2\delta_2\frac{d^2 u_e^\la}{d\la^2}+\la\delta_3\frac{du_e^\la}{d\la}+\delta_4 u_e^\la\big]\frac{d u_e^\la}{d\la}\\
=&\frac{d}{d\la}\int_{\R^{n+1}_+\cap\pa B_1}\theta_1^b\big[-5\frac{d^3 u_e^\la}{d\la^3}\frac{d u_e^\la}{d\la}+(20-5\delta_1)\la^3\frac{d^2 u_e^\la}{d\la^2}\frac{d u_e^\la}{d\la}\big]\\
&+\int_{\R^{n+1}_+\cap\pa B_1}\theta_1^b\big[(5\delta_1-20)\la^3(\frac{d^2 u_e^\la}{d\la^2})^2-5\delta_3\la(\frac{d u_e^\la}{d\la})^2\big]\\
&+\int_{\R^{n+1}_+\cap\pa B_1}\theta_1^b\big[5\la^4\frac{d^3 u_e^\la}{d\la^3}\frac{d^2 u_e^\la}{d\la^2}+(15\delta_1-60-5\delta_2)\frac{d^2 u_e^\la}{d\la^2}\frac{d u_e^\la}{d\la}-5\delta_4\frac{d u_e^\la}{d\la}u_e^\la\big].
\endaligned\ee
Here we have use that
\be\nonumber\aligned
-\la^4 f''''f'=\big[-5\la^4f'''f'+20\la^3f''f'\big]'-20\la^3(f'')^2-60\la^2f''f'+5\la^4f'''f''
\endaligned\ee
and
\be\nonumber\aligned
-\la^3 f'''f'=\big[-\la^3f''f'\big]'+3\la^2 f''f'+\la^3 (f'')^2.
\endaligned\ee
Now we add up $I_1,I_2, I_3$ and further integrate by part, we can get the term $\mathcal{I}$.
\be\label{8I}\aligned
\mathcal{I}:=& I_1+I_2+I_3\\
=&\frac{d}{d\la}\int_{\R^{n+1}_+\cap\pa B_1}\theta_1^b\big[\la^5\frac{d^4 u_e^\la}{d\la^4}\frac{d u_e^\la}{d\la}-2\la^5\frac{d^3 u_e^\la}{d\la^3}\frac{d^2 u_e^\la}{d\la^2}\\
&+(\delta_1-6)\la^4\frac{d^3 u_e^\la}{d\la^3}\frac{d u_e^\la}{d\la}+(24-6\delta_1+\delta_2)\la^3\frac{d^2 u_e^\la}{d\la^2}\frac{d u_e^\la}{d\la}\\
&+(9\delta_1-3\delta_2-36)\la^2\frac{d u_e^\la}{d\la}\frac{d u_e^\la}{d\la}\\
&+(8-\delta_1)\la^4(\frac{d^2 u_e^\la}{d\la^2})^2-\delta_4\la\frac{d u_e^\la}{d\la}u_e^\la-2\delta_4(u_e^\la)^2\big]\\
&+\int_{\R^{n+1}_+\cap\pa B_1}\theta_1^b\Big(
2\la^5(\frac{d^3u_e^\la}{d\la^3})^2+(10\delta_1-2\delta_2-56)\la^3(\frac{d^2u_e^\la}{d\la^2})^2\\
&+(-18\delta_1+\delta_2-4\delta_3+2\delta_4+72)\la(\frac{du_e^\la}{d\la})^2
\Big).
\endaligned\ee
Since $u_e^\la(X)=\la^{\frac{2s}{p-1}}u_e(\la X)$, we have the following
\be\nonumber\aligned
\la^4&\frac{d^4 u_e^\la}{d\la^4}=\la^{\frac{2s}{p-1}}\big[\frac{2s}{p-1}(\frac{2s}{p-1}-1)(\frac{2s}{p-1}-2)(\frac{2s}{p-1}-3)u_e(\la X)\\
&+\frac{8s}{p-1}(\frac{2s}{p-1}-1)(\frac{2s}{p-1}-2)r\la\pa_r u_e(\la X)\\
&+\frac{12s}{p-1}(\frac{2s}{p-1}-1)r^2\la^2\pa_{rr}u_e(\la X)\\
&+\frac{8s}{p-1}r^3\la^3\pa_{rrr}u_e(\la X)+r^4\la^4\pa_{rrrr}u_e(\la X)\big],\\
\endaligned\ee
and
\be\nonumber\aligned
\la^3\frac{d^3 u_e^\la}{d\la^3}&=\la^{\frac{2s}{p-1}}\big[\frac{2s}{p-1}(\frac{2s}{p-1}-1)(\frac{2s}{p-1}-2)u_e(\la X)\\
&+\frac{6s}{p-1}(\frac{2s}{p-1}-1)r\la\pa_r u_e(\la X)\\
&+\frac{6s}{p-1}r^2\la^2 \pa_{rr} u_e(\la X)+r^3\la^3\pa_{rrr} u_e(\la X)\big],
\endaligned\ee
\be\nonumber\aligned
&\la^2\frac{d^2 u_e^\la}{d\la^2}\\
&=\la^{\frac{2s}{p-1}}\big[\frac{2s}{p-1}(\frac{2s}{p-1}-1) u_e(\la X)+\frac{4s}{p-1}r\la\pa_r u_e(\la X)+r^2\la^2\pa_{rr}u_e(\la X)\big]
\endaligned\ee
and
\be\nonumber\aligned
\la\frac{d u_e^\la}{d\la}=\la^{\frac{2s}{p-1}}\big[\frac{2s}{p-1}u_e(\la X)+r\la\pa_r u_e(\la X)\big].
\endaligned\ee
Hence,  by scaling we have the following derivatives:
\be\nonumber\aligned
\frac{d}{d\la}&\int_{\R^{n+1}_+\cap\pa B_1}\theta_1^b\la^5\frac{d^4 u_e^\la}{d\la^4}\frac{d u_e^\la}{d\la}\\
&=\frac{d}{d\la}\int_{\R^{n+1}_+\cap\pa B_\la}\la^{2s\frac{p+1}{p-1}-n-5}y^b\big[\frac{2s}{p-1}(\frac{2s}{p-1}-1)(\frac{2s}{p-1}-2)(\frac{2s}{p-1}-3)u_e\\
&\quad+\frac{8s}{p-1}(\frac{2s}{p-1}-1)(\frac{2s}{p-1}-2)\la\pa_r u_e+\frac{12s}{p-1}(\frac{2s}{p-1}-1)\la^2\pa_{rr}u_e\\
&\quad+\frac{8s}{p-1}\la^3\pa_{rrr}u_e+\la^4\pa_{rrrr}u_e\big]\big[\frac{2s}{p-1}u_e+r\la\pa_r u_e\big];
\endaligned\ee

\be\label{8Iscaling}\aligned
\frac{d}{d\la}&\int_{\R^{n+1}_+\cap\pa B_1}\theta_1^b\la^5\frac{d^3 u_e^\la}{d\la^3}\frac{d^2 u_e^\la}{d\la^2}\quad\quad\quad\quad\\
&=\frac{d}{d\la}\int_{\R^{n+1}_+\cap\pa B_\la}\la^{2s\frac{p+1}{p-1}-n-5}y^b\big[\frac{2s}{p-1}(\frac{2s}{p-1}-1)(\frac{2s}{p-1}-2)u_e\\
&+\frac{6s}{p-1}(\frac{2s}{p-1}-1)\la\pa_r u_e
+\frac{6s}{p-1}\la^2 \pa_{rr} u_e+\la^3\pa_{rrr} u_e\big]\\
&\big[\frac{2s}{p-1}(\frac{2s}{p-1}-1) u_e+\frac{4s}{p-1}\la\pa_r u_e+\la^2\pa_{rr}u_e\big];
\endaligned\ee

\be\nonumber\aligned
\frac{d}{d\la}&\int_{\R^{n+1}_+\cap\pa B_1}\theta_1^b\la^4\frac{d^3 u_e^\la}{d\la^3}\frac{d u_e^\la}{d\la}\quad\quad\quad\quad\\
&=\frac{d}{d\la}\int_{\R^{n+1}_+\cap\pa B_\la}\la^{2s\frac{p+1}{p-1}-n-5}y^b\big[\frac{2s}{p-1}(\frac{2s}{p-1}-1)(\frac{2s}{p-1}-2)u_e\\
&+\frac{6s}{p-1}(\frac{2s}{p-1}-1)\la\pa_r u_e
+\frac{6s}{p-1}\la^2 \pa_{rr} u_e+\la^3\pa_{rrr} u_e\big]\\
&\big[\frac{2s}{p-1}u_e+\la\pa_r u_e\big];
\endaligned\ee

\be\nonumber\aligned
\frac{d}{d\la}&\int_{\R^{n+1}_+\cap\pa B_1}\theta_1^b\la^3\frac{d^2 u_e^\la}{d\la^2}\frac{d u_e^\la}{d\la}\quad\quad\quad\quad\\
&=\frac{d}{d\la}\int_{\R^{n+1}_+\cap\pa B_\la}\la^{2s\frac{p+1}{p-1}-n-5}y^b
\big[\frac{2s}{p-1}(\frac{2s}{p-1}-1) u_e\\
&+\frac{4s}{p-1}\la\pa_r u_e
+\la^2\pa_{rr}u_e][\frac{2s}{p-1}u_e+\la\pa_r u_e\big]
\endaligned\ee
and
\be\nonumber\aligned
\frac{d}{d\la}&\int_{\R^{n+1}_+\cap\pa B_1}\theta_1^b\la^2\frac{d u_e^\la}{d\la}\frac{d u_e^\la}{d\la}\quad\quad\quad\quad\\
&=\frac{d}{d\la}\int_{\R^{n+1}_+\cap\pa B_\la}\la^{2s\frac{p+1}{p-1}-n-5}y^b
\big[\frac{2s}{p-1}u_e+\la\pa_r u_e\big]^2.
\endaligned\ee
Further,
\be\nonumber\aligned
\frac{d}{d\la}&\int_{\R^{n+1}_+\cap\pa B_1}\theta_1^b\la^4\frac{d^2 u_e^\la}{d\la^2}\frac{d^2 u_e^\la}{d\la^2}\\
&=\frac{d}{d\la}\int_{\R^{n+1}_+\cap\pa B_\la}\la^{2s\frac{p+1}{p-1}-n-5}y^b
\big[\frac{2s}{p-1}(\frac{2s}{p-1}-1) u_e\\
&+\frac{4s}{p-1}\la\pa_r u_e+\la^2\pa_{rr}u_e\big]^2,
\endaligned\ee

\be\nonumber\aligned
\frac{d}{d\la}&\int_{\R^{n+1}_+\cap\pa B_1}\theta_1^b\la\frac{d u_e^\la}{d\la}u_e^\la\\
&=\frac{d}{d\la}\int_{\R^{n+1}_+\cap\pa B_\la}\la^{2s\frac{p+1}{p-1}-n-5}y^b
\big[\frac{2s}{p-1}u_e+\la\pa_r u_e\big]u_e,
\endaligned\ee
and
\be\nonumber\aligned
\frac{d}{d\la}&\int_{\R^{n+1}_+\cap\pa B_1}\theta_1^b u_e^\la=\frac{d}{d\la}\int_{\R^{n+1}_+\cap\pa B_\la}\la^{2s\frac{p+1}{p-1}-n-5}y^b
u_e^2.
\endaligned\ee


\subsection{The computations of $J_i,K_i,L_i (i=1,2,3)$ and $\mathcal{J},\mathcal{K},\mathcal{L}$}
Firstly,
\be\label{8j1}\aligned
J_1:=&\int_{\R^{n+1}_+\cap\pa B_1}\la\theta_1^b \frac{d}{d\la}J(u_e^\la)\frac{du_e^\la}{d\la}=\int_{\R^{n+1}_+\cap\pa B_1}\la\theta_1^b J(\frac{du_e^\la}{d\la})\frac{du_e^\la}{d\la}\\
=&\la\int_{\R^{n+1}_+\cap\pa B_1} \mathbf{div}_{S^n}\big(\theta_1^b\nabla_{S^n}(\theta_1^{-b}\mathbf{div}_{S^n}(\theta_1^b\nabla_{S^n}\frac{du_e^\la}{d\la}))\big)\frac{du_e^\la}{d\la}\\
=&-\la\int_{\R^{n+1}_+\cap\pa B_1}\theta_1^b\nabla_{S^n}\big(\theta_1^{-b} \mathbf{div}_{S^n}(\theta_1^b \nabla_{S^n}\frac{du_e^\la}{d\la})\big)\nabla_{S^n}\frac{du_e^\la}{d\la}\\
=&\la\int_{\R^{n+1}_+\cap\pa B_1}\theta_1^{-b}\big[\mathbf{div}_{S^n}(\theta_1^b\nabla_{S^n}\frac{du_e^\la}{d\la})\big]^2\\
=&\la\int_{\R^{n+1}_+\cap\pa B_1}\theta_1^b(\Delta_{S^n}\frac{du_e^\la}{d\la})^2,
 \endaligned\ee
here we have used integrate by part formula on the unit  sphere $S^n$.
\be\label{8j2}\aligned
J_2:=&-\la\int_{\R^{n+1}_+\cap\pa B_1}\theta_1^b J(u_e^\la)\frac{d^2 u_e^\la}{d\la^2}\\
=&-\la\int_{\R^{n+1}_+\cap\pa B_1}\mathbf{div}_{S^n}\big(\theta_1^b\nabla_{S^n}(\theta_1^{-b}\mathbf{div}_{S^n}(\theta_1^b\nabla_{S^n}u_e^\la))\big)\frac{d^2 u_e^\la}{d\la^2}\\
=&\la\int_{\R^{n+1}_+\cap\pa B_1}\theta_1^b\nabla_{S^n}(\theta_1^{-b}\mathbf{ div}_{S^n}(\theta_1^b \nabla_{S^n}u_e^\la)\nabla_{S^n}\frac{d^2u_e^\la}{d\la^2}\\
=&-\la\int_{\R^{n+1}_+\cap\pa B_1} \theta_1^{-b}\mathbf{div}_{S^n}(\theta_1^b \nabla_{S^n}u_e^\la) \frac{d^2}{d\la^2}\mathbf{div}_{S^n}(\theta_1^b\nabla_{S^n}u_e^\la)\\
=&\frac{d}{d\la}\int_{\R^{n+1}_+\cap\pa B_1}-\la \theta_1^{-b}\big[\mathbf{div}_{S^n}(\theta_1^b \nabla_{S^n}u_e^\la)\big]\frac{d}{d\la}\big[\mathbf{div}_{S^n}(\theta_1^b \nabla_{S^n}u_e^\la)\big]\\
 &\quad+\int_{\R^{n+1}_+\cap\pa B_1}\theta_1^{-b}\mathbf{ div}_{S^n}(\theta_1^b \nabla_{S^n}u_e^\la)\cdot \frac{d}{d\la}\mathbf{div}_{S^n}(\theta_1^b \nabla_{S^n}u_e^\la)\\
 &\quad+\la\int_{\R^{n+1}_+\cap\pa B_1} \theta_1^{-b}\big[\frac{d}{d\la}\mathbf{div}_{S^n}(\theta_1^b \nabla_{S^n}u_e^\la)\big]^2\\
=&\frac{d}{d\la}\Big(\int_{\R^{n+1}_+\cap\pa B_1}\theta_1^b\Big(
-\frac{1}{2}\la\frac{d}{d\la}(\Delta_{S^n}u_e^\la)^2+\frac{1}{2}(\Delta_{S^n}u_e^\la)^2
\Big)\Big)\\
&+\int_{\R^{n+1}_+\cap\pa B_1}\theta_1^{b}\la(\Delta_{S^n}\frac{du_e^\la}{d\la})^2,
 \endaligned\ee
here we denote that $g= \mathbf{div}_{S^n}(\theta_1^b \nabla_{S^n}u_e^\la),g'=\frac{d}{d\la}\mathbf{div}_{S^n}(\theta_1^b \nabla_{S^n}u_e^\la)$
and we have used that
\be\nonumber
-\la g g'=\big[-g g'\big]'+g g'+\la (g')^2=\big[-g g'+\frac{1}{2}g^2\big]'+\la(g')^2.
\ee
Further,
\be\label{8j3}\aligned
J_3:=&-5\int_{\R^{n+1}_+\cap\pa B_1}\theta_1^b J(u_e^\la)\frac{du_e^\la}{d\la}\\
=&-5\int_{\R^{n+1}_+\cap\pa B_1}\mathbf{div}_{S^n}\big(\theta_1^b \nabla_{S^n}(\theta_1^{-b}\mathbf{div}_{S^n}(\theta_1^b\nabla_{S^n}u))\big)\frac{du_e^\la}{d\la}\\
=&5\int_{\R^{n+1}_+\cap\pa B_1}\theta_1^b \nabla_{S^n}\big(\theta_1^{-b}\mathbf{div}_{S^n}(\theta_1^b\nabla_{S^n}u_e^\la)\big)\nabla_{S^n}\frac{d u_e^\la}{d\la}\\
=&-5\int_{\R^{n+1}_+\cap\pa B_1}\theta_1^{-b}\mathbf{ div}_{S^n}(\theta_1^b\nabla_{S^n}u_e^\la)\frac{d}{d\la}\mathbf{div}_{S^n}(\theta_1^b\nabla_{S^n}u_e^\la)\\
=&-\frac{5}{2}\frac{d}{d\la}\int_{\R^{n+1}_+\cap\pa B_1}\theta_1^b(\Delta_{S^n}u_e^\la)^2.
 \endaligned\ee
Therefore, combine with \eqref{8j1}, \eqref{8j2} and \eqref{8j3}, we get that
\be\label{8j}\aligned
\mathcal{J}:=&J_1+J_2+J_3\\
=&2\la\int_{\R^{n+1}_+\cap\pa B_1}\theta_1^{-b}\big[\frac{d}{d\la}\mathbf{div}_{S^n}(\theta_1^b\nabla_{S^n}u_e^\la)\big]^2\\
 &\quad-4\int_{\R^{n+1}_+\cap\pa B_1}\theta_1^{-b}\mathbf{div}_{S^n}(\theta_1^b\nabla_{S^n}u)\frac{d}{d\la}\mathbf{div}_{S^n}(\theta_1^b\nabla_{S^n}u)\\
 &\quad+\frac{d}{d\la}\int_{\R^{n+1}_+\cap\pa B_1}-\la \theta_1^b\big[\mathbf{div}_{S^n}(\theta_1^b \nabla_{S^n}u_e^\la)\big]\frac{d}{d\la}\big[\mathbf{div}_{S^n}(\theta_1^b \nabla_{S^n}u_e^\la)\big]\\
=&2\la\int_{\R^{n+1}_+\cap\pa B_1}\theta_1^{-b}\big[\frac{d}{d\la}\mathbf{div}_{S^n}(\theta_1^b\nabla_{S^n}u_e^\la)\big]^2\\
 &\quad-2\frac{d}{d\la}\int_{\R^{n+1}_+\cap\pa B_1}\theta_1^{-b}\mathbf{div}_{S^n}(\theta_1^b\nabla_{S^n}u_e^\la)\mathbf{div}_{S^n}(\theta_1^b\nabla_{S^n}u_e^\la)\\
 &\quad+\frac{d}{d\la}\int_{\R^{n+1}_+\cap\pa B_1}-\la \theta_1^{-b}\big[\mathbf{div}_{S^n}(\theta_1^b \nabla_{S^n}u_e^\la)\big]\frac{d}{d\la}\big[\mathbf{div}_{S^n}(\theta_1^b \nabla_{S^n}u_e^\la)\big]\\
=&\frac{d}{d\la}\int_{\R^{n+1}_+\cap\pa B_1}\theta_1^{b}\Big(-2(\Delta_{S^n}u_e^\la)^2-\frac{1}{2}\la\frac{d}{d\la}(\Delta_{S^n}u_e^\la)^2\Big)\\
 &\quad+2\la\int_{\R^{n+1}_+\cap\pa B_1}\theta_1^{b}(\Delta_{S^n}u_e^\la)^2.
\endaligned\ee
Note that
\be\label{8Jscaling}\aligned
&\frac{d}{d\la}\int_{\R^{n+1}_+\cap\pa B_1}\theta_1^b \big[\theta_1^{-b}\mathbf{div}_{S^n}(\theta_1^b\nabla_{S^n}u_e^\la)\big]^2\\
&=\frac{d}{d\la}\int_{\R^{n+1}_+\cap\pa B_\la}\la^{2s\frac{p+1}{p-1}-n-5}\big(\la^2\Delta_b u_e-\la^2\pa_{rr}u_e-(n+b)\la\pa_r u_e\big)^2,
\endaligned\ee
and
\be\nonumber\aligned
&\frac{d}{d\la}\int_{\R^{n+1}_+\cap\pa B_1}\theta_1^b\la \frac{d}{d\la}\big[\theta_1^{-b}\mathbf{div}_{S^n}(\theta_1^b\nabla_{S^n}u_e^\la)\big]^2\\
&=\frac{d}{d\la}\int_{\R^{n+1}_+\cap\pa B_\la}\la^{2s\frac{p+1}{p-1}-n-4}\frac{d}{d\la}\big(\la^2\Delta_b u_e-\la^2\pa_{rr}u_e-(n+b)\la\pa_r u_e\big)^2.
\endaligned\ee
Next we compute $K_1,K_2,K_3$ and $\mathcal{K}$.
\be\label{8k1}\aligned
&K_1\\
:=&\la\int_{\R^{n+1}_+\cap\pa B_1}\theta_1^b\frac{d}{d\la}K(u_e^\la)\frac{du_e^\la}{d\la}\\
=&2\la\int_{\R^{n+1}_+\cap\pa B_1}\mathbf{div}_{S^n}\big(\theta_1^b\nabla_{S^n}(\frac{d}{d\la}(\pa_{rr}+(n+b-2)\pa_r)u_e^\la)\big)\frac{du_e^\la}{d\la}\\
=&2\la\int_{\R^{n+1}_+\cap\pa B_1}\mathbf{div}_{S^n}\big(\theta_1^b\nabla_{S^n}(\la^3\frac{d^3 u_e^\la}{d\la^3}+\la^2(\alpha+2)\frac{d^2 u_e^\la}{d\la^2}+\la(\alpha+\beta)\frac{d u_e^\la}{d\la})\big)\frac{du_e^\la}{d\la}\\
=&-2\la\int_{\R^{n+1}_+\cap\pa B_1}\theta_1^b\nabla_{S^n}\big(\la^3\frac{d^3 u_e^\la}{d\la^3}+\la^2(\alpha+2)\frac{d^2 u_e^\la}{d\la^2}+\la(\alpha+\beta)\frac{d u_e^\la}{d\la}\big)\nabla_{S^n}\frac{du_e^\la}{d\la}\\
=&\frac{d}{d\la}\int_{\R^{n+1}_+\cap\pa B_1}-\la^3\theta_1^b\frac{d}{d\la}\big(\frac{d}{d\la}\nabla_{S^n} u_e^\la\big)^2
+(2-2\alpha)\la^2\int_{\R^{n+1}_+\cap\pa B_1}\theta_1^b (\frac{d}{d\la}\nabla_{S^n} u_e^\la)\\
\cdot &(\frac{d^2}{d\la^2}\nabla_{S^n} u_e^\la)
+2\la^3\int_{\R^{n+1}_+\cap\pa B_1}\theta_1^b(\frac{d^2}{d\la^2}\nabla_{S^n} u_e^\la)^2\\
&-(2\alpha+2\beta)\la\int_{\R^{n+1}_+\cap\pa B_1}\theta_1^b(\frac{d}{d\la}\nabla_{S^n} u_e^\la)^2.
\endaligned\ee
Here we denote that $h=\nabla_{S^n} u_e^\la,h'=\frac{d}{d\la}\nabla_{S^n} u_e^\la$, and have used that
\be\nonumber
-\la^3h'h'''=\big[-\frac{\la^3}{2}\frac{d}{d\la}(h')^2\big]'+3\la^2 h'h''+\la^3(h'')^2.
\ee
Next,
\be\label{8k2}\aligned
K_2:=&-\la\int_{\R^{n+1}_+\cap\pa B_1}\theta_1^b K(u_e^\la)\frac{d^2 u_e^\la}{d\la^2}\\
=&-2\la\int_{\R^{n+1}_+\cap\pa B_1}\mathbf{div}_{S^n}(\theta_1^b\nabla_{S^n}\big(\la^2\frac{d^2 u_e^\la}{d\la^2}+\la\alpha\frac{d u_e^\la}{d\la}+\beta u_e^\la)\big)\frac{d^2 u_e^\la}{d\la^2}\\
=&2\la\int_{\R^{n+1}_+\cap\pa B_1}\theta_1^b\nabla_{S^n}\big(\la^2\frac{d^2 u_e^\la}{d\la^2}+\la\alpha\frac{d u_e^\la}{d\la}+\beta u_e^\la\big)\nabla_{S^n}\frac{d^2 u_e^\la}{d\la^2}\\
=&\frac{d}{d\la}\int_{\R^{n+1}_+\cap\pa B_1}\theta_1^b\big[2\beta\la\nabla_{S^n}u_e^\la \frac{d}{d\la}\nabla_{S^n}u_e^\la-\beta(\nabla_{S^n}u_e^\la)^2\big]\\
&+2\la^3\int_{\R^{n+1}_+\cap\pa B_1}\theta_1^b(\frac{d^2}{d\la^2}\nabla_{S^n}u_e^\la)^2-2\la\beta\int_{\R^{n+1}_+\cap\pa B_1}\theta_1^b(\frac{d}{d\la}\nabla_{S^n}u_e^\la)^2\\
&+2\la^2\alpha\int_{\R^{n+1}_+\cap\pa B_1}\theta_1^b\frac{d}{d\la}\nabla_{S^n}u_e^\la\frac{d^2}{d\la^2}\nabla_{S^n}u_e^\la.
\endaligned\ee
Here we have used that
\be\nonumber
2\la h h''=[2\la h h'-h^2]'-2\la(h')^2.
\ee
Further,
\be\label{8k3}\aligned
K_3:=&-5\int_{\R^{n+1}_+\cap\pa B_1}\theta_1^b K(u_e^\la)\frac{du_e^\la}{d\la}\\
=&-10\int_{\R^{n+1}_+\cap\pa B_1}\mathbf{div}_{S^n}\big(\theta_1^b\nabla_{S^n}(\la^2\frac{d^2 u_e^\la}{d\la^2}+\la\alpha\frac{du_e^\la}{d\la}+\beta u_e^\la)\big)\frac{du_e^\la}{d\la}\\
=&10\int_{\R^{n+1}_+\cap\pa B_1}\theta_1^b\nabla_{S^n}\big(\la^2\frac{d^2 u_e^\la}{d\la^2}+\la\alpha\frac{du_e^\la}{d\la}+\beta u_e^\la\big)\nabla_{S^n}\frac{du_e^\la}{d\la}\\
=&\frac{d}{d\la}\big[5\beta\int_{\R^{n+1}_+\cap\pa B_1}\theta_1^b \nabla_{S^n} u_e^\la \nabla_{S^n} u_e^\la\big]+10\la\alpha\int_{\R^{n+1}_+\cap\pa B_1}\theta_1^b\big(\frac{d}{d\la}\nabla_{S^n} u_e^\la\big)^2\\
&+10\la^2\int_{\R^{n+1}_+\cap\pa B_1}\theta_1^b \frac{d}{d\la}\nabla_{S^n} u_e^\la \frac{d^2}{d\la^2}\nabla_{S^n} u_e^\la.
\endaligned\ee
Now combine with \eqref{8k1}, \eqref{8k2} and \eqref{8k3}, we get that
\be\label{8k}\aligned
\mathcal{K}:=&K_1+K_2+K_3\\
=&\frac{d}{d\la}\int_{\R^{n+1}_+\cap\pa B_1}\theta_1^b\Big[-\la^3\frac{d}{d\la}(\frac{d}{d\la}\nabla_{S^n} u_e^\la)^2\\
&+2\beta\la\nabla_{S^n} u_e^\la\frac{d}{d\la}\nabla_{S^n} u_e^\la+4\beta(\nabla_{S^n} u_e^\la)^2+6\la^2(\nabla_{S^n}\frac{du_e^\la}{d\la})^2\Big]\\
&+4\la^3\int_{\R^{n+1}_+\cap\pa B_1}\theta_1^b(\frac{d^2}{d\la^2}\nabla_{S^n} u_e^\la)^2\\
&+(8\alpha-4\beta-12)\la\int_{\R^{n+1}_+\cap\pa B_1}\theta_1^b(\frac{d}{d\la}\nabla_{S^n} u_e^\la)^2.\\
\endaligned\ee
Notice that by scaling we have
\be\label{8Kscaling}\aligned
&\frac{d}{d\la}\int_{\R^{n+1}_+\cap\pa B_1}\theta_1^b(\nabla_{S^n} u_e^\la)^2\\
&=\frac{d}{d\la}\int_{\R^{n+1}_+\cap\pa B_\la}\la^{2s\frac{p+1}{p-1}-n-5}y^b\big[\la^2|\nabla u_e|^2-\la^2|\pa_r u_e|^2\big].
\endaligned\ee
\be\nonumber\aligned
&\frac{d}{d\la}\int_{\R^{n+1}_+\cap\pa B_1}\theta_1^b\la\frac{d}{d\la}(\nabla_{S^n} u_e^\la)^2\\
&=\frac{d}{d\la}\int_{\R^{n+1}_+\cap\pa B_\la}\la^{2s\frac{p+1}{p-1}-n-4}y^b\frac{d}{d\la}\big[\la^2|\nabla u_e|^2-\la^2|\pa_r u_e|^2\big]
\endaligned\ee
and
\be\nonumber\aligned
&\frac{d}{d\la}\int_{\R^{n+1}_+\cap\pa B_1}\theta_1^b\la^3\frac{d}{d\la}(\frac{d}{d\la}\nabla_{S^n} u_e^\la)^2\\
&=\frac{d}{d\la}\int_{\R^{n+1}_+\cap\pa B_\la}\la^{2s\frac{p+1}{p-1}-n-4}y^b\frac{d}{d\la}\big[\frac{2s}{p-1}\la\nabla u_e+\la^2\nabla \pa_r u_e\big]^2.
\endaligned\ee
Finally, we compute $\mathcal{L}$.
\be\nonumber\aligned
&L_1:=\int_{\R^{n+1}_+\cap\pa B_1}\la\theta_1^b\frac{d}{d\la}L(u_e^\la)\frac{du_e^\la}{d\la}\\
&=-2(n+b-3)\la\int_{\R^{n+1}_+\cap\pa B_1}\mathbf{div}_{S^n}(\theta_1^b\nabla_{S^n}\frac{du_e^\la}{d\la})\frac{du_e^\la}{d\la}\\
&=2(n+b-3)\la\int_{\R^{n+1}_+\cap\pa B_1}\theta_1^b(\nabla_{S^n}\frac{du_e^\la}{d\la})^2;
\endaligned\ee

\be\nonumber\aligned
&L_2:=\int_{\R^{n+1}_+\cap\pa B_1}-\la\theta_1^bL(u_e^\la)\frac{d^2u_e^\la}{d\la^2}\\
&=2(n+b-3)\la\int_{\R^{n+1}_+\cap\pa B_1}\mathbf{div}_{S^n}(\theta_1^b\nabla_{S^n}u_e^\la)\frac{d^2 u_e^\la}{d\la^2}\\
&=-2(n+b-3)\int_{\R^{n+1}_+\cap\pa B_1}\la\theta_1^b\nabla_{S^n}u_e^\la\frac{d^2}{d\la^2}\nabla_{S^n}u_e^\la\\
&=-(n+b-3)\int_{\R^{n+1}_+\cap\pa B_1}\theta_1^b\frac{d}{d\la}\big[2\la\nabla_{S^n}u_e^\la\nabla_{S^n}\frac{du_e^\la}{d\la}-(\nabla_{S^n}u_e^\la)^2\big]\\
&+2(n+b-3)\la\int_{\R^{n+1}_+\cap\pa B_1}\theta_1^b|\nabla_{S^n}\frac{du_e^\la}{d\la}|^2;\\
\endaligned\ee

\be\nonumber\aligned
&L_3:=\int_{\R^{n+1}_+\cap\pa B_1}-5\theta_1^bL(u_e^\la)\frac{du_e^\la}{d\la}\\
&=10(n+b-3)\int_{\R^{n+1}_+\cap\pa B_1}\mathbf{div}_{S^n}(\theta_1^b\nabla_{S^n}u_e^\la)\frac{du_e^\la}{d\la}\\
&=-10(n+b-3)\int_{\R^{n+1}_+\cap\pa B_1}\theta_1^b\nabla_{S^n}u_e^\la\nabla_{S^n}\frac{du_e^\la}{d\la}\\
&=-5(n+b-3)\frac{d}{d\la}\int_{\R^{n+1}_+\cap\pa B_1}\theta_1^b[\nabla_{S^n}u_e^\la]^2.\\
\endaligned\ee
Hence,
\be\nonumber\aligned
&\mathcal{L}:=L_1+L_2+L_3\\
&=-(n+b-3)\frac{d}{d\la}\int_{\R^{n+1}_+\cap\pa B_1}\theta_1^b\big[\la\frac{d}{d\la}(\nabla_{S^n}u_e^\la)^2\big]\\
&\quad-4(n+b-3)\frac{d}{d\la}\int_{\R^{n+1}_+\cap\pa B_1}\theta_1^b[\nabla_{S^n}u_e^\la]^2\\
&\quad+4(n+b-4)\la\int_{\R^{n+1}_+\cap\pa B_1}\theta_1^b|\nabla_{S^n}\frac{du_e^\la}{d\la}|^2\\
&=-(n+b-3)\frac{d}{d\la}\int_{\R^{n+1}_+\cap\pa B_\la}\la^{2s\frac{p+1}{p-1}-n-4}y^b\frac{d}{d\la}\big[\la^2|\nabla u_e|^2-\la^2|\pa_r u_e|^2\big]\\
&\quad-4(n+b-3)\frac{d}{d\la}\int_{\R^{n+1}_+\cap\pa B_\la}\la^{2s\frac{p+1}{p-1}-n-5}y^b\big[\la^2|\nabla u_e|^2-\la^2|\pa_r u_e|^2\big].
\endaligned\ee
By rescaling, we have
\be\label{8Lscaling}\aligned
&\frac{d}{d\la}\int_{\R^{n+1}_+\cap\pa B_1}\theta_1^b\big[\la\frac{d}{d\la}(\nabla_{S^n}u_e^\la)^2\big]\\
&=\frac{d}{d\la}\int_{\R^{n+1}_+\cap\pa B_\la}\la^{2s\frac{p+1}{p-1}-n-4}y^b\frac{d}{d\la}\big[\la^2|\nabla u_e|^2-\la^2|\pa_r u_e|^2\big];\\
&\frac{d}{d\la}\int_{\R^{n+1}_+\cap\pa B_1}\theta_1^b[\nabla_{S^n}u_e^\la]^2\\
&=\frac{d}{d\la}\int_{\R^{n+1}_+\cap\pa B_\la}\la^{2s\frac{p+1}{p-1}-n-5}y^b\big[\la^2|\nabla u_e|^2-\la^2|\pa_r u_e|^2\big].
\endaligned\ee

\subsection{The term $\overline{E}_{d_1}$}
Notice that on the boundary $\pa B_1$,
\be\nonumber\aligned
&v_e^\la=\Delta_b u_e^\la\\
&=\la^2\frac{d^2u_e^\la}{d\la^2}+(n+b-\frac{4s}{p-1})\la\frac{du_e^\la}{d\la}+\frac{2s}{p-1}(1+\frac{2s}{p-1}-n-b)u_e^\la+\Delta_\theta u_e^\la\\
&:=\la^2\frac{d^2u_e^\la}{d\la^2}+\alpha_0\la\frac{du_e^\la}{d\la}+\beta_0 u_e^\la+\Delta_\theta u_e^\la.
\endaligned\ee
Integrate by part, it follows that

\be\label{8ed1ed1}\aligned
&\int_{\pa B_1}y^b\la(\frac{dv_e^\la}{d\la})^2=\int_{\pa B_1}\Big(\la^5(\frac{d^3u_e^\la}{d\la^3})^2+(\alpha_0^2-2\alpha_0-2\beta_0-4)\la^3(\frac{d^2u_e^\la}{d\la^2})^2\\
&\quad\quad\quad\quad\quad\quad\quad\quad+(-\alpha_0^2+\beta_0^2+2\alpha_0+2\beta_0)\la(\frac{du_e^\la}{d\la})^2\Big)\\
&+\int_{\pa B_1}\Big(-2\la^3(\nabla_\theta \frac{d^2u_e^\la}{d\la^2})^2+(10-2\beta_0)\la(\nabla_\theta \frac{du_e^\la}{d\la})^2\Big)
+\int_{\pa B_1}\la(\Delta_\theta \frac{du_e^\la}{d\la})^2\\
&+\frac{d}{d\la}\Big(\int_{\pa B_1}\sum_{0\leq i,j\leq 2,i+j\leq2}{c^1_{i,j}}\la^{i+j}\frac{d^iu_e^\la}{d\la^i}\frac{d^ju_e^\la}{d\la^j}
+\sum_{0\leq s,t\leq 2,s+t\leq2}{c^2_{s,t}}\la^{s+t}\frac{d^su_e^\la}{d\la^s}\frac{d^tu_e^\la}{d\la^t}\Big),
\endaligned\ee
where $c^i_{i,j},c^2_{s,t}$ depending on   $a,b$ hence on  $p,n$.

\vskip0.2in

\noindent{\bf Proof of Theorem \ref{8monoid}}. Notice that the equation \eqref{8overlineE}, combine with the estimates on $\mathcal{I},\mathcal{J},\mathcal{K},\mathcal{L}$ and \eqref{8ed1ed1}, we obtain Theorem \ref{8monoid}.\hfill $\Box$

\section{Energy estimates and Blow down analysis }

 In this section, we do some energy estimates for the solutions of \eqref{8LE}, which are important when we perform a
 blow-down analysis in the next section.

 \subsection{Energy estimates}

\bl Let $u$ be a solution of \eqref{8LE} and $u_e$ satisfy \eqref{8LEE}, then there exists a positive constant $C$ such that
\be\label{8estimate1}\aligned
&\int_{\pa\R^{n+1}_+}|u_e|^{p+1}\eta^6+\int_{\R^{n+1}_+}y^b|\nabla\Delta_b u_e|^2\eta^6\\
&\leq C\big[\int_{\R^{n+1}_+}y^b|\Delta_b u_e|^2\eta^4|\nabla\eta|^2+\int_{\R^{n+1}_+}y^b|\nabla u_e|^2 \frac{|\Delta_b\eta^6|^2}{\eta^6}
+\int_{\R^{n+1}_+}y^bu_e^2\frac{|\nabla\Delta_b \eta^6|^2}{\eta^6}\\
&\quad +\int_{\R^{n+1}_+}y^b|\nabla u_e|^2\eta^2|\nabla \eta|^4
+\int_{\R^{n+1}_+}y^b|\nabla^2 u_e|^2\eta^4|\nabla\eta|^2\\
&\quad +\int_{\R^{n+1}_+}y^b|\nabla u_e|^2\eta^4|\nabla^2\eta|^2\big].
\endaligned\ee
\el
\bp
Multiply the equation \eqref{8LEE} with $y^b u_e\eta^6$, where $\eta$ is a test function, we get that
\be\aligned
0=&\int_{\R^{n+1}_+}y^b u_e\eta^6\Delta_b^3 u_e=\int_{\R^{n+1}_+} u_e\eta^6 \mathbf{div}(y^b\nabla\Delta_b^2u_e)\\
=&-\int_{\pa\R^{n+1}_+} u_e\eta^6 \frac{\pa}{\pa y}\Delta_b^2 u_e-\int_{\R^{n+1}_+}y^b\nabla(u_e\eta^6)\nabla\Delta_b^2u_e\\
=&C_{n,s}\int_{\pa\R^{n+1}_+}|u_e|^{p+1}\eta^6-\int_{\pa\R^{n+1}_+}y^b\frac{\pa (u_e\eta^6)}{\pa y}\Delta_b^2 u_e+\int_{\R^{n+1}_+}y^b\Delta_b(u_e\eta^6)\Delta_b^2u_e\\
=&C_{n,s}\int_{\pa\R^{n+1}_+}|u_e|^{p+1}\eta^6+\int_{\R^{n+1}_+}y^b\Delta_b(u_e\eta^6)\Delta_b^2u_e\\
=&C_{n,s}\int_{\pa\R^{n+1}_+}|u_e|^{p+1}\eta^6-\int_{\pa\R^{n+1}_+}\Delta_b(u_e\eta^6)y^b\frac{\pa \Delta_b u_e}{\pa y}\\
&-\int_{\R^{n+1}_+}y^b\nabla(\Delta_b(u_e\eta^6))\nabla\Delta_b u_e\\
=&C_{n,s}\int_{\pa\R^{n+1}_+}|u_e|^{p+1}\eta^6-\int_{\R^{n+1}_+}y^b\nabla(\Delta_b(u_e\eta^6))\nabla\Delta_b u_e.\\
\endaligned\ee
Hence, we have
\be\label{8eee}
C_{n,s}\int_{\pa\R^{n+1}_+}|u_e|^{p+1}\eta^6=\int_{\R^{n+1}_+}y^b\nabla(\Delta_b(u_e\eta^6))\nabla\Delta_b u_e.
\ee
Since $\Delta_b(\xi\eta)=\eta\Delta_b\xi+\xi\Delta_b\eta+2\nabla\xi\nabla\eta$, we have
\be\nonumber\aligned
\Delta_b(u_e\eta^6)=\eta^6\Delta_b u_e+u_e\Delta_b\eta^6+12\eta^5\nabla u_e\nabla\eta,
\endaligned\ee
therefore,
\be\label{8zhankai}\aligned
\nabla\Delta_b(u_e\eta^6)\nabla\Delta_b u_e=&6\eta^5\Delta_b u_e\nabla\eta\nabla\Delta_b u_e+(\eta)^6(\nabla\Delta_b u_e)^2+\Delta_b\eta^6\nabla u_e\nabla\Delta_b u_e\\
&+u_e\nabla\Delta_b \eta^6\nabla\Delta_b u_e+60\eta^4(\nabla\eta\nabla\Delta_b u_e)(\nabla u_e\nabla\eta)\\
&+12\eta^5\sum_{i,j}\pa_{ij}u_e\pa_i\eta\pa_j \Delta_b u_e+12\eta^5\sum_{i,j}\pa_i u_e\pa_{ij}\eta\pa_j\Delta_b u_e.
\endaligned\ee
here $\pa_j(j=1, ..., n, n+1)$ denote the derivatives  with respect to  $x_1,...,x_n,y$ respectively.
A similar way can be applied to deal with the following term $|\nabla\Delta_b (u_e\eta^3)|^2$.
On the other hand, by the stability condition, we have
\be\label{8stable1}\aligned
p\int_{\R^n}|u|^{p+1}\eta^6\leq\int_{\R^n}\int_{\R^n}\frac{\Big(u(x)\eta^3(x)-u(y)\eta^3(y)\Big)^2}{|x-y|^{n+2s}}=
\frac{1}{C_{n,s}}\int_{\R^{n+1}_+}y^b|\nabla\Delta_b(u_e\eta^3)|^2.
\endaligned\ee
(Here we notice that $u_e(x,0)=u(x)$, see Theorem \ref{ththYang}, \eqref{Yangmore})\\
Combine with \eqref{8eee}, \eqref{8zhankai} and \eqref{8stable1}, we have
\be\nonumber\aligned
&\int_{\R^{n+1}_+}y^b|\nabla\Delta_b u_e|^2\eta^6\\
&\leq  C\varepsilon\int_{\R^{n+1}_+}y^b(\nabla\Delta_b u_e)^2\eta^6+C(\varepsilon)\big[\int_{\R^{n+1}_+}y^b|\Delta_b u_e|^2\eta^4|\nabla\eta|^2\\
&\quad +\int_{\R^{n+1}_+}y^b|\nabla u_e|^2 (\frac{|\Delta_b\eta^6|^2}{\eta^6}+\eta^4|\nabla^2\eta|^2)\\
&\quad +\int_{\R^{n+1}_+}y^bu_e^2\frac{|\nabla\Delta_b \eta^6|^2}{\eta^6}
+\int_{\R^{n+1}_+}y^b|\nabla u_e|^2\eta^2|\nabla \eta|^4
+\int_{\R^{n+1}_+}y^b|\nabla^2 u_e|^2\eta^4|\nabla\eta|^2\big],\\
\endaligned\ee
we can select $\varepsilon$ so small   that $C\varepsilon\leq\frac{1}{2}$. Combine with \eqref{8eee} and \eqref{8zhankai}, we obtain our conclusion.
\ep


\bc\label{8hhhh}
 Let $u$ be a solution of \eqref{8LE} and $u_e$ satisfy \eqref{8LEE}, then

\be\nonumber\aligned
&\int_{\pa\R^{n+1}_+\cap B_{R/2}}|u_e|^{p+1}+\int_{\R^{n+1}_+\cap B_{R/2}}y^b(\nabla\Delta_b u_e)^2\\
&\leq C\big[R^{-6}\int_{\R^{n+1}_+\cap B_R}y^b u_e^2+R^{-4}\int_{\R^{n+1}_+\cap B_R}y^b|\nabla u_e|^2\\
&+R^{-2}\int_{\R^{n+1}_+\cap B_R}y^b(|\Delta_b u_e|^2+|\nabla^2 u_e|^2)\big].
\endaligned\ee
\ec
\bp
We let $\eta=\xi^m$ where $m>1$ in the estimate \eqref{8estimate1}. We have
\be\nonumber\aligned
&\int_{\pa\R^{n+1}_+}|u_e|^{p+1}\xi^{6m}+\int_{\R^{n+1}_+}y^b|\nabla\Delta_b u_e|^2\xi^{6m}\\
&\leq C[\int_{\R^{n+1}_+}y^b(|\Delta_b u_e|^2+|\nabla^2 u_e|^2)\xi^{6m-2}|\nabla\xi|^2\\
&\quad +\int_{\R^{n+1}_+}y^b|\nabla u_e|^2\xi^{6m-4}(|\nabla^2\xi|^2+|\nabla\xi|^4)+\int_{\R^{n+1}_+}y^b u_e^2\xi^{6m-6}|\nabla^3\xi|^2].
\endaligned\ee
Let $\xi=1$ in $B_{R/2}$ and $\xi=0$ in $B_R^C$, satisfying $|\nabla \xi |\leq\frac{C}{R}$, then we have the desired estimates.
\ep

\vskip0.1in
\bl\label{8lemma3}
Suppose that $u$ is a solution of \eqref{8LE} which is stable outside some ball $B_{R_0}\subset\R^n$. For $\eta\in C_c^\infty(\R^n\backslash\overline{B_{R_0}})$ and $x\in\R^n$, define
\be
\rho(x)=\int_{\R^n}\frac{(\eta(x)-\eta(y))^2}{|x-y|^{n+2s}}dy.
\ee
Then
\be
\int_{\R^n}|u|^{p+1}\eta^2dx+\int_{\R^n}\int_{\R^n}\frac{|u(x)\eta(x)-u(y)\eta(y)|^2}{|x-y|^{n+2s}}dxdy\leq C\int_{\R^n} u^2\rho dx.
\ee
\el

\bl\label{8lemma4}
Let $m>n/2$ and $x\in\R^n$. Set
\be\label{8420}
\rho(x):=\int_{\R^n}\frac{(\eta(x)-\eta(y))^2}{|x-y|^{n+2s}}dy\;\;\hbox{where}\;\;\eta(x)=(1+|x|^2)^{-m/2}.
\ee
Then there is a constant $C=C(n,s,m)>0$ such that
\be C^{-1}(1+|x|^2)^{-n/2-s}\leq \rho(x)\leq C(1+|x|^2)^{-n/2-s}.
\ee
\el

\bc\label{8corr1}
Suppose that $m>n/2$, $\eta$ is given by \eqref{8420} and $R>R_0>1$. Define
\be\label{8rhor}
\rho_R(x)=\int_{\R^n}\frac{(\eta_R(x)-\eta_R(y))^2}{|x-y|^{n+2s}}dy,\;\;\hbox{where}\;\;\eta_R(x)=\eta(x/R)\psi(x/R)
\ee
for a standard test function $\psi$ that $\psi\in C^{\infty}(\R^n)$, $0\leq\psi\leq1,\psi=0$ on $B_1$ and $\psi=1$ on $\R^n\setminus B_2$.
Then there exists a constant $C>0$ such that
\be\nonumber
\rho_R(x)\leq C\eta^2(x/R)|x|^{-(n+2s)}+R^{-2s}\rho(x/R).
\ee
\ec

\bl\label{8lemma5}
Suppose that $u$ is a   solution of \eqref{8LE} which  is stable outside a ball $B_{R_0}$. Consider $\rho_R$  which  is defined in \eqref{8rhor} for
$n/2<m<n/2+s(p+1)/2$. Then there exists a constant $C>0$ such that
\be\nonumber
\int_{\R^n}u^2\rho_R\leq C(\int_{B_{3R_0}}u^2\rho_R+R^{n-2s\frac{p+1}{p-1}})
\ee
for any $R>3R_0$.
\el
The proofs of Lemma \ref{8lemma3}, Corollary \ref{8corr1}, Lemma \ref{8lemma4} and Lemma \ref{8lemma5} are similar to that of  Lemmas $2.1, 2.2, 2.4$ in \cite{Wei0=1} and we omit the details here.

\bl\label{8estimate111}
Suppose that $p\neq\frac{n+2s}{n-2s}$. Let $u$ be a solution of \eqref{8LE} which is stable outside a ball $B_{R_0}$ and $u_e$ satisfy \eqref{8LEE}. Then
there exists a constant $C>0$ such that
\be\nonumber
\int_{B_R}y^b u_e^2\leq CR^{n+6-2s\frac{p+1}{p-1}},\quad \int_{B_R}y^b |\nabla u_e|^2\leq CR^{n+4-2s\frac{p+1}{p-1}},
\ee
\be\nonumber
\int_{B_R}y^b (|\nabla^2 u_e|^2+|\Delta_b u_e|^2)\leq CR^{n+2-2s\frac{p+1}{p-1}}.
\ee
\el
\bp
Recall that the Possion formula for the fractional equation for the case $0<s<1$  (see \cite{Caffarelli2007}), we can generalize the
 expression formula to the general case with non-integer positive real number. Therefore,
\be\nonumber
{u}_e(x,y)=C_{n,s}\int_{\R^n}u(z)\frac{y^{2s}}{(|x-z|^2+y^2)^{\frac{n+2s}{2}}}dz.
\ee
Then we have
\be\label{8pa0}
|{u}_e(x,y)|^2\leq C\int_{\R^n}u^2(z)\frac{y^{2s}}{(|x-z|^2+y^2)^{\frac{n+2s}{2}}}dz,
\ee
and
\be\nonumber
\pa_y{u}_e(x,y)=C_{n,s}\int_{\R^n}u(z)\big[\frac{2sy^{2s-1}}{(|x-z|^2+y^2)^{\frac{n+2s}{2}}}
-\frac{(n+2s)y^{2s+1}}{(|x-z|^2+y^2)^{\frac{n+2s+2}{2}}}\big]dz,
\ee
also
\be\nonumber
\pa_{x_j}{u}_e(x,y)=-C_{n,s}\int_{\R^n}u(z)\frac{(n+2s)(x_j-z_j)y^{2s}}{(|x-z|^2+y^2)^{\frac{n+2s+2}{2}}}dz,
\ee
for $j=1,2,...,n$.
Hence by H\"older's inequality we have
\be\label{8pa1}
|\nabla u_e(x,y)|^2\leq C\int_{\R^n}\frac{u^2(z)y^{2s-2}}{(|x-z|^2+y^2)^{\frac{n+2s}{2}}}dz.
\ee
By a straightforward calculation we have
\be\nonumber\aligned
\pa_{x_jx_j}{u}_e(x,y)=&C_{n,s}\int_{\R^n}u(z)\big[\frac{(n+2s)(n+2s+2)(x_j-z_j)^2y^{2s}}{(|x-z|^2+y^2)^{\frac{n+2s+4}{2}}}\\
&-\frac{(n+2s)y^{2s}}{(|x-z|^2+y^2)^{\frac{n+2s+2}{2}}}\big]dz,
\endaligned\ee

\be\nonumber\aligned
\pa_{x_jy}{u}_e(x,y)=&C_{n,s}\int_{\R^n}u(z)\big[\frac{(n+2s)(n+2s+2)(x_j-z_j)^2y^{2s+1}}{(|x-z|^2+y^2)^{\frac{n+2s+4}{2}}}\\
&-\frac{2s(n+2s)(x_j-z_j)^2y^{2s-1}}{(|x-z|^2+y^2)^{\frac{n+2s+2}{2}}}\big],
\endaligned\ee
and
\be\nonumber\aligned
\pa_{yy}{u}_e(x,y)=&C_{n,s}\int_{\R^n}u(z)\big[\frac{2s(2s-1)y^{2s-2}}{(|x-z|^2+y^2)^{\frac{n+2s}{2}}}\\
&-\frac{(n+2s)(4s+1)y^{2s}}{(|x-z|^2+y^2)^{\frac{n+2s+2}{2}}}
+\frac{(n+2s)(n+2s+2)y^{2s+2}}{(|x-z|^2+y^2)^{\frac{n+2s+4}{2}}}\big].
\endaligned\ee
Therefore, we have
\be\nonumber
|\nabla^2 u_e(x,y)|+|\Delta_b u_e(x,y)|\leq C\int_{\R^n}|u(z)|\frac{y^{2s-2}}{(|x-z|^2+y^2)^{\frac{n+2s}{2}}}dz.
\ee
Hence,
\be\label{8pa2}
|\nabla^2 u_e(x,y)|^2+|\Delta_b u_e(x,y)|^2\leq C\int_{\R^n}u^2(z)\frac{y^{2s-4}}{(|x-z|^2+y^2)^{\frac{n+2s}{2}}}dz.
\ee
Now we turn to estimate the following integration, which provides a unify way to deal with our desired estimates.

Define
\be\label{8ak}\aligned
A_k:=&\int_{|x|\leq R,z\in\R^n}u^2(z)\big[\int_0^R\frac{y^{2k+1}}{(|x-z|^2+y^2)^{\frac{n+2s}{2}}}dy\big]dzdx\\
=&\int_{|x|\leq R,z\in\R^n}u^2(z)\big[\int_0^{R^2}\frac{\alpha^k}{(|x-z|^2+\alpha)^{\frac{n+2s}{2}}}d\alpha]dzdx\\
\leq&\frac{1}{2}\int_{|x|\leq R,z\in\R^n}u^2(z)[\int_0^{R^2}\frac{d\alpha}{(|x-z|^2+\alpha)^{\frac{n+2s}{2}}-k}\big]dzdx\\
=&\frac{1}{2}(\frac{n+2s}{2}-k)\int_{|x|\leq R,z\in\R^n}u^2(z)\big[(|x-z|^2)^{k-\frac{n+2s}{2}+1}\\
&-(|x-z|^2+R^2)^{k-\frac{n+2s}{2}+1}\big],\\
\endaligned\ee
where $k=0,1,2$.
Split the integral to $|x-z|\leq 2R$ and $|x-z|>2R$, for the case of $|x-z|\leq 2R$, we see that

\be\nonumber\aligned
&\int_{|x|\leq R,|x-z|\leq 2R}u^2(z)\Big[(|x-z|^2)^{k-\frac{n+2s}{2}+1}-(|x-z|^2+R^2)^{k-\frac{n+2s}{2}+1}\Big]\\
&\leq\int_{|x|\leq R,|x-z|\leq 2R}u^2(z)\Big[(|x-z|^2)^{k-\frac{n+2s}{2}+1}\Big]\\
&\leq C R^{2k-2s+2}\int_{|z|\leq 3R}u^2(z)dz\\
&\leq R^{2k-2s+2}\Big(\int_{B_{3R}}|u|^{p+1}\eta_R^2\Big)^{2/(p+1)}\Big(\int_{B_{3R}}\eta_R^{-4/(p-1)} \Big)^{(p-1)/(p+1)}\\
&\leq C R^{2k-2s+2}\Big(\int_{B_{3R}}u^2(z)\rho_R(z)\Big)^{2/(p+1)}\\
&\leq C R^{n+2k+2-2s\frac{p+1}{p-1}}.
\endaligned
\ee
Here we have used Lemma \ref{8lemma3} and \ref{8lemma5}. For the case of $|x-z|>2R$, by the mean value theorem, we have
\be\nonumber\aligned
&\int_{|x|\leq R,|x-z|> 2R}u^2(z)\Big[(|x-z|^2)^{k-\frac{n+2s}{2}+1}-(|x-z|^2+R^2)^{k-\frac{n+2s}{2}+1}\Big]\\
&\leq R^2\int_{|x|\leq R,|x-z|> 2R}u^2(z)\Big[(|x-z|^2)^{k-\frac{n+2s}{2}}\Big]\\
&\leq C R^{n+2}\int_{|z|\geq R} u_e^2(z)|z|^{2k-n-2s}dz\\
&\leq C R^{n+2}\Big[\int_{|z|\geq R} (u^{p+1}_e(z))\Big]^{2/(p+1)}\Big(\int_{|z|\geq R}|z|^{(2k-n-2s)\frac{p+1}{p-1}}\Big)^{(p-1)/(p+1)}\\
&\leq C R^{n+2k+2-2s\frac{p+1}{p-1}},
\endaligned
\ee
here we have used Lemma \ref{8lemma3}.
Hence, we obtain that
\be\label{8akak}
A_k\leq C R^{n+2k+2-2s\frac{p+1}{p-1}},
\ee
where $C=C(n,s,p)$ independent of $R$.
Now, combine with \eqref{8pa0}, \eqref{8pa1} and \eqref{8pa2}, recall that $b=5-2s$, we have

\be\nonumber\aligned
&\int_{B_R}y^b u_e^2dxdy\leq A_2, \quad \int_{B_R}y^b|\nabla u_e|^2dxdy\leq A_1,\\
 &\int_{B_R}y^b\big(|\nabla^2 u_e|^2+|\Delta_b u_e|^2\big)dxdy\leq A_0.
\endaligned
\ee
Apply \eqref{8akak}, we finish our proof.
\ep

Combine Corollary \ref{8hhhh} and Lemma \ref{8lemma5}, we have the following lemma.
\bl\label{8lemmaE}
Let $u$ be a solution of \eqref{8LE} which is stable outside a ball $B_{R_0}$ and $u_e$ satisfy \eqref{8LEE}. Then there exists a positive constant $C$
such that
\be\nonumber\aligned
&\int_{\pa\R^{n+1}_+\cap B_R}|u_e|^{p+1}+R^{-6}\int_{\R^{n+1}_+\cap B_R}y^b |u_e|^2
+R^{-4}\int_{\R^{n+1}_+\cap B_R}y^b |\nabla u_e|^2\\
&+R^{-2}\int_{\R^{n+1}_+\cap B_R}y^b\big(|\Delta_b u_e|^2+|\nabla^2 u_e|^2\big)
+\int_{\R^{n+1}_+\cap B_R}y^b|\nabla\Delta_b u_e|^2
\leq C R^{n-2s\frac{p+1}{p-1}}.
\endaligned\ee
\el

\newpage
\subsection{ Blow down analysis and  the proof of Theorem \ref{8Liouvillec}}

  {\bf The proof of Theorem \ref{8Liouvillec}.}
 Suppose that $u$ is a solution of \eqref{8LE} which is stable outside the ball of radius $R_0$ and suppose that $u_e$ satisfies \eqref{8LEE}.
 In the subcritical case, i.e., $1<p<p_s(n)$, Lemma \ref{8lemma3} implies that $u\in \dot{H}^s(\R^n)\cap L^{p+1}(\R^n)$.
Multiplying \eqref{8LE} with $u$ and integrate, we obtain that
\be\label{841}
\int_{\R^n}|u|^{p+1}=\|u\|^2_{\dot{H}^s(\R^n)}.
\ee
Multiplying \eqref{8LE} with $u^\la(x)=u(\la x)$ yields
\be\nonumber
\int_{\R^n}|u|^{p-1} u^\la=\int_{\R^n} (-\Delta)^{s/2}u (-\Delta)^{s/2}u^\la=\la^s\int_{\R^n}w w_\la,
\ee
where $w=(-\Delta)^{s/2}u$. Following the ideas provided in \cite{Xavier2015,Xavier2014} and using the change of variable $z=\sqrt{\la}x$, we can get
the following Pohozaev identity
 \be\nonumber\aligned
-\frac{n}{p+1}\int_{\R^n}|u|^{p+1}&=\frac{2s-n}{2}\int_{\R^n}|w|^{2}
+\frac{d}{d\la}\int_{\R^n}w^{\sqrt{\la}}w^{1/\sqrt{\la}}dz \Big|_{\la=1}\\
&=\frac{2s-n}{2}\|u\|^2_{\dot{H}^s(\R^n)}.
\endaligned\ee
Hence, we have the following Pohozaev identity
 \be\nonumber\aligned
\frac{n}{p+1}\int_{\R^n}|u|^{p+1}=\frac{n-2s}{2}\|u\|^2_{\dot{H}^s(\R^n)}.
\endaligned\ee

For $p<p_s(n)$, this equality above together with \eqref{841} proves that $u\equiv0$.
For $p=p_s(n)$, this equality above means that the energy is finite. Further,  since $u\in \dot{H}^s(\R^n)$, apply the stability inequality with
test function $\psi=u\eta^2(\frac{x}{R})$, and let $R\rightarrow+\infty$ (where $\eta$ is cutoff function), then we get that
 \be\nonumber\aligned
p\int_{\R^n}|u|^{p+1}\leq\|u\|^2_{\dot{H}^s(\R^n)}.
\endaligned\ee
This together with \eqref{841} gives that $u\equiv0$.

\vskip0.1in

Now we consider the supercritical case, i.e., $p>\frac{n+2s}{n-2s}$, we perform the proof via a few steps.

\vskip0.1in

\noindent{\bf Step 1.}  $\lim_{\la\rightarrow\infty} E(u_e,0,\la)<\infty$.\\

From Theorem \ref{8Monotone} we know that $E$ is nondecreasing w.r.t. $\la$, so we only need to show that $E(u_e,0,\la)$ is bounded. Note that
\be\nonumber
E(u_e,0,\la)\leq\frac{1}{\la}\int_{\la}^{2\la} E(u_e,0,t)dt\leq\frac{1}{\la^2}\int_\la^{2\la}\int_{t}^{t+\la}E(u_e,0,\gamma)d\gamma dt.
\ee
From Lemma \ref{8lemmaE}, we have that
\be\nonumber\aligned
\frac{1}{\la^2}&\int_\la^{2\la}\int_{t}^{t+\la}\gamma^{2s\frac{p+1}{p-1}-n}\big[\int_{\R^{n+1}_+\cap B_\gamma}\frac{1}{2}y^b|\nabla\Delta_b u_e|^2dydx\\
&\quad\quad\quad\quad\quad\quad -\frac{C_{n,s}}{p+1}\int_{\pa\R^{n+1}_+\cap B_\gamma}|u_e|^{p+1}dx\big]d\gamma dt\leq C,
\endaligned\ee
where $C>0$ is independent of $\gamma$.
\be\label{8type1}\aligned
\frac{1}{\la^2}&\int_\la^{2\la}\int_{t}^{t+\la}\int_{\R^{n+1}_+\cap\pa B_\gamma}\gamma^{2s\frac{p+1}{p-1}-n-5}y^b\big[\frac{2s}{p-1}(\frac{2s}{p-1}-1)(\frac{2s}{p-1}-2)u_e\\
&+\frac{6s}{p-1}(\frac{2s}{p-1}-1)\gamma\pa_r u_e
+\frac{6s}{p-1}\gamma^2 \pa_{rr} u_e+\gamma^3\pa_{rrr} u_e\big]\\
&\quad\quad\big[\frac{2s}{p-1}(\frac{2s}{p-1}-1) u_e+\frac{4s}{p-1}\gamma\pa_r u_e+\gamma^2\pa_{rr}u_e\big]\\
&\leq C\frac{1}{\la^2}\int_\la^{2\la}\int_{t}^{t+\la}t^{2s\frac{p+1}{p-1}-n-5}\int_{\R^{n+1}_+\cap\pa B_\gamma}\\
&\quad y^b\big[u_e^2+\gamma^2(\pa_r u_e)^2+\gamma^4(\pa_{rr}u_e)^2+\gamma^6(\pa_{rrr}u_e)^2\big]\\
&\leq C\frac{1}{\la^2}\int_\la^{2\la}t^{2s\frac{p+1}{p-1}-n-5}\int_{\R^{n+1}_+\cap B_{3\la}}\\
&\quad y^b\big[u_e^2+\gamma^2(\pa_r u_e)^2+\gamma^4(\pa_{rr}u_e)^2+\gamma^6(\pa_{rrr}u_e)^2\big]\\
&\leq C \la^{n-2s\frac{p+1}{p-1}+6}\frac{1}{\la^2}\int_\la^{2\la}t^{2s\frac{p+1}{p-1}-n-5}dt\\
&\leq C
\endaligned\ee
and

\be\label{8type2}\aligned
&\Big|\frac{1}{\la^2}\int_\la^{2\la}\int_{t}^{t+\la}\int_{\R^{n+1}_+\cap\pa B_\gamma}\gamma^{2s\frac{p+1}{p-1}-n-4}y^b\frac{d}{d\gamma}(\gamma^2\Delta_b u_e-\gamma^2\pa_{rr}u_e-(n+b)\gamma\pa_r u_e)^2\Big|\\
&\leq\frac{1}{\la^2}\int_\la^{2\la}t^{2s\frac{p+1}{p-1}-n-5}\int_{t}^{t+\la}\int_{\R^{n+1}_+\cap\pa B_\gamma}y^b\big[2\gamma^2\Delta_b u_e-2\gamma^2\pa_{rr}u_e-(n+b)\gamma\pa_r u_e\big]\\
&\quad \big[\gamma^2\Delta_b u_e-\gamma^2\pa_{rr}u_e-(n+b)\gamma\pa_r u_e\big]\\
&\leq\frac{1}{\la^2}\int_\la^{2\la}t^{2s\frac{p+1}{p-1}-n-5}\int_{t}^{t+\la}\int_{\R^{n+1}_+\cap B_{3\la}}y^b[2\gamma^2\Delta_b u_e-2\gamma^2\pa_{rr}u_e-(n+b)\gamma\pa_r u_e]\\
&\quad\big[\gamma^2\Delta_b u_e-\gamma^2\pa_{rr}u_e-(n+b)\gamma\pa_r u_e\big]\\
&\leq C\frac{1}{\la^2}\int_\la^{2\la}t^{2s\frac{p+1}{p-1}-n-5}\int_{t}^{t+\la}\int_{\R^{n+1}_+\cap B_{3\la}}y^b
\big[\gamma^4(\Delta_b u_e)^2+\gamma^4(\pa_{rr}u_e)^2+\gamma^2\pa_r u_e\big]\\
&\leq C \la^{n-2s\frac{p+1}{p-1}+6}\frac{1}{\la^2}\int_\la^{2\la}t^{2s\frac{p+1}{p-1}-n-5}dt\\
&\leq C.
\endaligned\ee
\vskip0.1in
Integrate by part,  by  the scaling identity of section 3, for example \eqref{8Iscaling}, \eqref{8Jscaling}, \eqref{8Kscaling} and \eqref{8Lscaling},
we can treat the remaining    terms  by a  similar way as the estimates \eqref{8type1} and \eqref{8type2}.

\vskip0.12in

\noindent{\bf Step 2.}  There exists a sequence $\la_i\rightarrow\infty$ such that $(u_e^{\la_i})$ converges weakly
to a function $u_e^\infty$  in $H^3_{loc}(\R^n;y^bdxdy)$,  this is a direct consequence of Lemma \ref{8lemmaE}.

\vskip0.12in

\noindent{\bf Step 3.}  {\bf The function $u_e^\infty$ is homogeneous.}
Due to the scaling invariance of $E $  (i.e., $E(u_e,0,R\la)=E(u_e^{\la},0,R)$ )
and the monotonicity formula, for any given $R_2>R_1>0$, we see that
\be\nonumber\aligned
0=&\lim_{i\rightarrow\infty}\big(E(u_e,0,R_2\la_i)-E(u_e,0,R_1\la_i)\big)\\
=&\lim_{i\rightarrow\infty}\big(E(u_e^{\la_i},0,R_2)-E(u_e^{\la_i},0,R_1)\big)\\
\geq&\liminf_{i\rightarrow\infty}\int_{(B_{R_2}\setminus B_{R_1})\cap \R^{n+1}_+}
y^b r^{2s\frac{p+1}{p-1}-n-6}\big(\frac{2s}{p-1}{u_e^{\la_i}}+r\frac{\pa u_e^{\la_i}}{\pa r}\big)^2dydx\\
 \geq&\int_{(B_{R_2}\setminus B_{R_1})\cap \R^{n+1}_+}
y^b r^{2s\frac{p+1}{p-1}-n-6}\big(\frac{2s}{p-1}{u_e^{\infty}}+r\frac{\pa u_e^{\infty}}{\pa r}\big)^2dydx.\\
\endaligned\ee
In the last inequality we have used the weak convergence of  the sequence $(u_e^{\la_i})$ to the function $u_e^{\infty}$ in $H^3_{loc}(\R^n;y^bdxdy)$. This implies that
\be\nonumber
\frac{2s}{p-1}\frac{u_e^{\infty}}{r}+\frac{\pa u_e^{\infty}}{\pa r}=0\;\;\hbox{a.e.}\;\;\hbox{in}\;\;\R^{n+1}_+.
\ee
Therefore,  $u_e^\infty$ is homogeneous.

\vskip0.12in

\noindent{\bf Step 4}. {  \bf $u_e^\infty=0$.}
This is a direct consequence of Theorem $3.1$ in \cite{Wei1=2}.

\vskip0.12in

\noindent{\bf Step 5}. $(u_e^{\la_i})$ converges strongly to zero in $H^3(B_R\setminus B_\varepsilon;y^bdxdy)$ and $(u_e^{\la_i})$ converges strongly in
$L^{p+1}(\pa\R^{n+1}_+\cap (B_R\setminus B_\varepsilon))$ for all $R>\varepsilon>0$. These are consequent results  of Lemma \ref{8lemmaE} and Theorem $1.5$ in \cite{Fabes1982}.

\vskip0.12in

\noindent{\bf Step 6}. $u_e=0$. Note that
\be\nonumber\aligned
\overline{E}(u_e,\la)=&\overline{E}(u_e^\la,1)\\
=&\frac{1}{2}\int_{\R^{n+1}_+\cap B_1}y^b|\nabla\Delta_b u_e^\la|^2dxdy-\frac{C_{n,s}}{p+1}\int_{\pa\R^{n+1}_+\cap B_1}|u_e^\la|^{p+1}dx\\
=&\frac{1}{2}\int_{\R^{n+1}_+\cap B_\varepsilon}y^b|\nabla\Delta_b u_e^\la|^2dxdy-\frac{C_{n,s}}{p+1}\int_{\pa\R^{n+1}_+\cap B_\varepsilon}|u_e^\la|^{p+1}dx\\
&+\frac{1}{2}\int_{\R^{n+1}_+\cap( B_1\setminus B_\varepsilon)}y^b|\nabla\Delta_b u_e^\la|^2dxdy-\frac{C_{n,s}}{p+1}\int_{\pa\R^{n+1}_+\cap( B_1\setminus B_\varepsilon)}|u_e^\la|^{p+1}dx\\
=&\varepsilon^{n-2s\frac{p+1}{p-1}}\overline{E}(u_e,\la\varepsilon)+\frac{1}{2}\int_{\R^{n+1}_+\cap( B_1\setminus B_\varepsilon)}y^b|\nabla\Delta_b u_e^\la|^2dxdy\\
&-\frac{C_{n,s}}{p+1}\int_{\pa\R^{n+1}_+\cap( B_1\setminus B_\varepsilon)}|u_e^\la|^{p+1}dx.\\
\endaligned\ee
Letting $\la\rightarrow+\infty$ and then $\varepsilon\rightarrow0$, we deduce that $\lim_{\la\rightarrow+\infty}\overline{E}(u_e,\la)=0$.
Using the monotonicity of $E$,
\be\aligned
E(u_e,\la)&\leq\frac{1}{\la}\int_{\la}^{2\la} E(t)dt\leq \sup_{[\la,2\la]} \overline{E}+C\frac{1}{\la}\int_\la^{2\la}[E-\overline{E}]\\
&\leq\sup_{[\la,2\la]} \overline{E}+C\frac{1}{\la} \int_\la^{2\la}\la^{2s\frac{p+1}{p-1}-n-5}\int_{\R^{n+1}_+\cap\pa B_\la}y^b
 \big[(u_e)^2+\la^2|\nabla u_e|^2\\
 &\quad +\la^4(|\Delta_b u_e|^2+|\nabla^2 u_e|^2)+\la^6|\nabla\Delta_b u_e|^2\big]\\
&=\sup_{[\la,2\la]} \overline{E}+C\frac{1}{\la}\la^{2s\frac{p+1}{p-1}-n-5}\int_{\R^{n+1}_+\cap( B_{2\la}\backslash B_\la)}y^b
\big[(u_e)^2+\la^2|\nabla u_e|^2\\
 &\quad +\la^4(|\Delta_b u_e|^2+|\nabla^2 u_e|^2)+\la^6|\nabla\Delta_b u_e|^2\big]\\
 &=\sup_{[\la,2\la]} \overline{E}+C\frac{1}{\la}\int_{\R^{n+1}_+\cap( B_{2}\backslash B_1)}y^b
\big[(u_e^\la)^2+|\nabla u_e^\la|^2\\
 &\quad +|\Delta_b u_e^\la|^2+|\nabla^2 u_e^\la|^2+|\nabla\Delta_b u_e^\la|^2\big]\\
\endaligned\ee
and so $\lim_{\la\rightarrow\infty} E(u_e,\la)=0$. Since $u$ is smooth, we also have $E(u_e,0)=0$. Since $E$ is monotone, $E\equiv0$
and so $\overline{u}$ must be homogenous, a contradiction unless $u_e\equiv0$.

\section{Algebraic analysis: The proof of Theorem \ref{8Monotonem}}


Let $k:=\frac{2s}{p-1}$ and $m:=n-2s$. By  a direct calculation,  we obtain  that
\be\label{8A1A2}\aligned
A_1&= -10k^2+10mk-m^2+12m+25,\\
A_2&=3k^4-6mk^3+(3m^2-12m-30)k^2+(12m^2+30m)k+9m^2+36m+27,\\
B_1&=-6k^2+6mk+12m+30.
\endaligned\ee

Notice that our supercritical condition $p>\frac{n+2s}{n-2s}$ is equivalent to $0<k<\frac{n-2s}{2}=\frac{m}{2}$.
Next, we have the following lemma which yields the sign of $A_2$ and $B_1$.

\bl\label{8A2-0} If $p>\frac{n+2s}{n-2s}$, then $A_2>0$ and $B_1>0$.
\el

\bp From  \eqref{8A1A2}, we derive that
\be\label{8A2}
A_2=3(k+1)(k+3)(k-(m+1))(k-(m+3)),
\ee
and the roots of $B_1=0$ are

\be\nonumber
\frac{1}{2}m-\frac{1}{2}\sqrt{m^2+8m+20},\quad \frac{1}{2}m+\frac{1}{2}\sqrt{m^2+8m+20}.
\ee
Recall that $p>\frac{n+2s}{n-2s}$ is equivalent to $0<k<\frac{m}{2}$, we get the conclusion.
\ep

 To show monotonicity formula,  we  proceed   to prove the following inequality. That is,  there exist real numbers $c_{i,j}$ and positive real number
 $\epsilon$ such that

\be\aligned
\label{8keyinequality}
&3\la^5(\frac{d^3u^\la}{d\la^3})^2+A_{1}\la^3\Big(\frac{d^2u^\la}{d\la^2}\Big)^2+A_{2}\la(\frac{du^\la}{d\la})^2\\
&\geq \epsilon\la(\frac{du^\la}{d\la})^2+\frac{d}{d\la}
\Big(\sum_{0\leq i,j\leq2}c_{i,j}\la^{i+j}\frac{d^iu^\la}{d\la^i}\frac{d^ju^\la}{d\la^j}
\Big).
\endaligned\ee
 To deal with the rest of the dimensions, we employ the second idea: we find nonnegative constants $d_1, d_2$ and constants $c_1, c_2$ such that we have the following Jordan form decomposition:
\be\label{8quad}\aligned
&3\la^5(f''')^2+A_1\la^3(f'')^2+A_2\la(f')^2=3\la(\la^2f'''+c_1\la f'')^2+d_1\la(\la f''+c_2 f')^2\\
&\;\quad\quad\quad\quad\quad\quad\quad\quad\quad\quad\quad\quad+d_2\la(f')^2
+\frac{d}{d\la}(\sum_{i,j}e_{i,j}\la^{i+j}f^{(i)}f^{(j)}),
\endaligned\ee
where the  unknown constants are to be determined.

\bl
\label{860011}
Let $ p>\frac{n+2s}{n-2s}$ and $A_1$ satisfy
\be
\label{8A11}
 A_1+12 >0,
 \ee
 then there exist  nonnegative numbers $d_1,d_2$,    and real numbers $c_1,c_2,e_{i,j}$ such that
the  differential inequality (\ref{8quad}) holds.

 \el

 \bp Since $$4\la^4f'''f''=\frac{d}{d\la}(2\la^4(f'')^2)-8\la^3(f'')^2$$ and
 $$2\la^2f''f'=\frac{d}{d\la}(\la^2(f')^2)-2\la(f')^2,$$ by comparing the coefficients of $\lambda^3(f'')^2$ and $\lambda (f')^2,$  we have that

\be\nonumber
d_1=A_1-3c_1^2+12c_1, \quad d_2=A_2-(c_2^2-2c_2)(A_1-3c_1^2+12c_1).
\ee
In particular, $$\max_{c_1}{d_1(c_1)}=A_1+12   \hbox{  and the critical point is }   c_1=2. $$
Since $A_2>0$,  we select that $c_1=2,c_2=0$.  Hence,  in this case,
by a direct calculation we see that $d_1=A_1+12>0$.  Then we get  the conclusion.
\ep

We conclude from Lemma \ref{860011} that if $A_1+12>0$ then  (\ref{8keyinequality}) holds.  This implies that when $m<6+\sqrt{73}, p>\frac{n+2s}{n-2s}$ or $m\geq 6+\sqrt{73}$ and
\be
\label{pee}
\frac{n+2s}{n-2s}<p<\frac{5m+20s-\sqrt{15m^2+120m+370}}{5m-\sqrt{15m^2+120m+370}},
\ee
 then (\ref{8keyinequality}) holds.
\vskip0.2in

Let
\be\label{727=1}
p_m(n):=\begin{cases}
+\infty\;\;\;\;\;\;\;\;&\hbox{if}\;  n<2s+6+\sqrt{73},\\
\frac{5n+10s-\sqrt{15(n-2s)^2+120(n-2s)+370}}{5n-10s-\sqrt{15(n-2s)^2+120(n-2s)+370}}&\hbox{if}\; n\geq2s+6+\sqrt{73}.\\
\end{cases}
\ee
Combining  all the lemmas of this section,  we obtain the Theorem \ref{8Monotonem}.

\vskip0.2in

\vskip0.16in

Now we proceed to prove Theorem \ref{8Monotone}. From  Corollary 1.1 of  \cite{LWZ}, we know that
 if $n>2s,s>0,p>\frac{n+2s}{n-2s}$,  then there exists  $n_0(s)$, where $\frac{1}{\sqrt{n}}<a_{n,s}<\frac{1}{2}\frac{n-2s}{\sqrt{n}}+\frac{1}{\sqrt{n}}$, such that   the inequality  \eqref{8gamma} always holds whenever
 $n\leq n_0(s)$; while when  $n> n_0(s)$, then  the inequality   \eqref{8gamma}  is true  if and only if
$$p<p_2:=\frac{n+2s-2-2a_{n,s}\sqrt{n}}{n-2s-2-2a_{n,s}\sqrt{n}},$$
where $n_0(s)$ is in fact the largest $n$ satisfying    $n-2s-2-2a_{n,s}\sqrt{n}\leq0$.
In particular, $    \frac{n+2s}{n-2s}<\frac{n+2s-4}{n-2s-4}<p_2<+\infty.$  Therefore, we introduce
\be\label{8pcn88}
p_c(n):=\begin{cases}
+\infty\;\;\;\;\;\;\;\;&\hbox{if}\;\;\;\;\;\;\;\; n\leq n_0(s),\\
\frac{n+2s-2-2a_{n,s}\sqrt{n}}{n-2s-2-2a_{n,s}\sqrt{n}}&\hbox{if}\;\;\;\;\;\;\;\;n> n_0(s).\\
\end{cases}
\ee

From \cite{LWZ=00}, we use the sharp estimate $n_0(s)<2s+8.998$ for $2<s<3$, then
\be\label{8n0s}
n_0(s)\leq 2s+8.998<2s+6+\sqrt{73}\simeq2s+14.544.
\ee
On the other hand, via the sharp estimate $a_{n,s}<1$ from \cite{LWZ=00}
\be\label{8compare}\aligned
\frac{5n+10s-\sqrt{15(n-2s)^2+120(n-2s)+370}}{5n-10s-\sqrt{15(n-2s)^2+120(n-2s)+370}}>\frac{n+2s-2-2a_{n,s}\sqrt{n}}{n-2s-2-2a_{n,s}\sqrt{n}}
\endaligned\ee
provided that   $s\in(2,3)$ and that
\be\label{728-mm}\aligned
225m^4-720m^3-17244m^2-29088m+7236>0,   \hbox{ where }  m=n-2s.
\endaligned\ee
The  \eqref{728-mm} holds  whenever  $m>11.12$, that is $n>2s+11.12$.
This combine with \eqref{8n0s} we obtain that $p_c(n)<p_m(n)$.
Therefore we get   Theorem \ref{8Monotone}. \hfill $\Box$

\newpage


\end{document}